\documentclass[11pt]{article}

\usepackage{amssymb,amsmath}
\usepackage{color}
\usepackage{enumerate}

\def\seq#1#2#3{#1_{#2},\,\ldots,#1_{#3}}

\def\N{{\mathbb N}}
\def\Z{{\mathbb Z}}
\def\Q{{\mathbb Q}}
\def\C{{\mathbb C}}
\def\P{{\mathbb P}}
\def\CP{{\mathbb CP}}
\def\k{\overline}
\def\w#1{\widetilde{#1}}

\def\h{\widehat}
\def\ord{\text{ord}}
\def\OO{{\mathcal O}}

\newtheorem{theo}{Theorem}
\newtheorem{lemm}{Lemma}
\newtheorem{coro}{Corollary}
\newtheorem{prop}{Proposition}

\newenvironment{defi}
{\smallskip\noindent{\bf Definition\/}:}{\smallskip\par}

\newenvironment{rema}
{\smallskip\noindent{\bf Remark\/}.}{\smallskip\par}

\title{Pencils and critical locus on normal surfaces}
\author{F. Delgado
\thanks{
Supported by the grant
MTM2012-36917-C03-01
with the help of FEDER Program.
The first author is  thankful to
the Institut Fourier, Universit\'e de
Grenoble I  for hospitality during the stages of this work.
\newline
Address: IMUVA (Instituto de Investigaci\'on en Matem\'aticas). Universidad de
Valladolid. Spain.
E-mail:
fdelgado\symbol{'100}agt.uva.es}
\and{H. Maugendre
\thanks{The second author is   thankful  to University of Valladolid
for hospitality during the stages of this work. \newline Address: Institut Fourier,
Universit\'{e} de Grenoble I, BP 74, F-38402 Saint-Martin
d'H\`{e}res, France.  E-mail: helene.maugendre\symbol{'100}ujf-grenoble.fr}
}}
\date{}

\begin{document}

\maketitle

\bigskip
\bigskip

\begin{abstract}
We study linear pencils of curves on normal
surface singularities. Using the minimal good resolution of the pencil, we describe
the topological type of  generic elements of the pencil and characterize the
behaviour of special elements. Then we show that the critical locus associated to
the pencil is  linked to  the special elements. This gives a decomposition of the
critical locus through the minimal good resolution and as a consequence,
information on the topological type of the critical locus.
\end{abstract}



\section{Introduction}

Let $ (Z,z)$ be a complex analytic normal surface, and let
 $\pi : (Z,z)\to (\C^2,0)$ be a  finite complex analytic morphism germ. We choose coordinates $(u,v)$ in
$(\C^2, 0)$ and denote $f:=u\circ \pi$ and $g:=v\circ \pi$.
We consider the meromorphic function $h:=f/g$ defined in a punctured neighbourhood $V$ of $z$ in $Z$.
It can be seen as a map $h: V\to \C\P^1$ defined by $h(x):=(f(x):g(x))$.
For  $w = (w_1:w_2)\in \C\P^1$,  the  closure of  $h^{-1}(w)$  defines the curve $ w_2f-w_1g=0$ on the surface $(Z,z)$.
The set $\Lambda:=\{ w_2f-w_1g,  w_1, w_2 \in \C\}$ is the {\it pencil} defined by $f$ and $g$.
We  denote $\phi_w$ the element  of the pencil $\Lambda$ equal to $w_2f-w_1g$. Its (non reduced) zero
locus, denoted by
$\Phi_w$,  is called the {\it  fibre} defined by $\phi_w$.

\bigskip
Such linear families of curves have been studied independently and
through different approach for $(Z,z)$ equal to $(\C^2,0)$  in  \cite{LW1}, \cite{DM} and
\cite{MM}. In the general case (it means $ (Z,z)$ a germ of normal complex analytic
surface which is not smooth anymore), L\^e D\~ung Tr\`ang and R. Bondil
give  in \cite{BL} a definition of general elements of the pencil which are
characterized by the minimality of their Milnor number. In \cite{B2} R. Bondil gives an
algebraic $\mu$-constant theorem for linear families of plane curves.
Other results have been obtained in the case where  $\pi$ is the restriction to
$(Z,z)$ of a linear projection of $(\C^n, 0)$ onto $(\C^2,0)$ (see \cite{B1}, \cite{BNP},
\cite{S}). At last, the topology of the morphism $\pi$ has been studied  in \cite{LMW}  and
\cite{Mi}.  In \cite{LMW}, the authors define rational quotients which are topological
invariants  of $(\pi, u,v)$  and give different ways to compute them. In \cite{Mi}, F.
Michel presents another proof of the topological invariance of this set of rational
numbers and moreover  she gives a decomposition of the critical locus of $\pi$  in
bunches linked to the set of invariants.

\bigskip
Let $\rho : (X,E)\to (Z,z)$ be a {\it good} resolution of  the singularity  $(Z,z)$. It
is a resolution of the singularity $(Z,z)$ such that the exceptional divisor is a
union of smooth projective curves with normal crossings.   In particular three irreducible components
of the exceptional divisor do not meet at the same point. The lifting $h\circ \rho$
is a  meromorphic function defined in a suitable neighbourhood of $E$ in $X$ but in
a finite set of points.

A {\it good resolution} $\rho$ of the pencil $ \Lambda $ is a good resolution of
the singularity  $(Z,z)$ in which  $h\circ \rho$ is a morphism and    the
exceptional divisor is a union of smooth projective curves with normal crossings. A good resolution of
the pencil $ \Lambda $ is said to be {\it minimal} if  and
only if by the contraction of any  rational component of self-intersection -1 of
the exceptional divisor
we do  not obtain a good resolution of $\Lambda$ anymore. We will see in section 2
that there  exists a unique minimal good resolution of $\Lambda$.

An irreducible component $E_\alpha$ of $E$ is called {\it dicritical} if the
restriction of $h\circ \rho$ to $E_\alpha$ is not constant.

\bigskip
Considering the minimal good resolution  $\rho : (Y,E)\to (Z,z)$  of the pencil
$\Lambda$, we define {\it special} and {\it generic} values of $\Lambda$ as follows.
Let us denote $\h h= h\circ \rho$ and $\cal D$ the union of the dicritical components of $E$.
We define the set of special zones $SZ(\Lambda)=\{ \Delta _i, i\in I\}$ where  $I $
is a finite subset of $\N$,  and  $\Delta_i$ is either a connected  components of
$\overline{E\setminus\cal D}$, either a critical point of the restriction of  $\h h$ to
$\cal D$, or an intersection point between two dicritical components.
Notice that $\h h_{|_{\Delta_i}}$ is constant.

\begin{defi} \rm \label{defi1}
The set of {\it special values} of $\Lambda$ is constituted of
 the values $\h h (\Delta_i)$ for $i\in I$.
A fibre associated to a special value  is called a {\it special fibre} of $\Lambda$.

The other values of $\C\P ^1$ are called {\it generic values} for the pencil $\Lambda$.
A  fibre associated to a generic value  is called a {\it generic fibre}  of $\Lambda$.
\end{defi}

We prove the  following results.

\begin{theo}\label{theoprincipal}
Let  $w,w'$ be   generic values for the pencil $\Lambda$, then
the fibers  $\Phi_w$ and $\Phi_{w'}$ have the same topological type.

Moreover, if $e\in \C\P^1$ is a special value for the pencil $\Lambda$, then the
fibers
$\Phi_w$ and $\Phi_e$ do not have the same topological type.
\end{theo}

The above definition and theorem generalize  some of the results  contained in \cite{LW1} (see theorem 4.1) where the authors study pencils defined on $\C^2$.
Going on studying the topology of the pencil we prove the following result which extend to the case of normal surfaces the second item of theorems 1, 2, 3 of \cite{DM}  which deals with pencils defined on $\C^2$.

\begin{theo}\label{specialfibre}
Let  $\rho$  be the minimal good resolution    of the pencil
$\Lambda$,  $\Delta\in SZ(\Lambda)$, and let $e\in \C  \P^1$. Then, the strict transform  of $\Phi_e$ by $\rho$ intersects
$\Delta$ if and only if $\Phi_e$ is special and
$\h{h}(\Delta)=e$.
\end{theo}

In a second part we are interested in understanding the
behaviour of the critical locus of the map $\pi$.
We denote  by $I_z(\ ,\ )$ the local intersection multiplicity at $z$ (see section 2.1).
 We prove the following result which generalize the third item of theorems 1, 2, 3 of \cite{DM}.

\begin{theo} \label{critloc}
Let  $\rho : (Y,E)\to (Z,z)$ be the minimal good resolution of
$\Lambda$.  For each element $\Delta\in SZ(\Lambda)$
there exists an irreducible
component  of   the critical locus $C(\pi )$  of $\pi $ such that its  strict
transform by $\rho$ intersects $\Delta $.

Moreover for each  branch $\Gamma$ of $C(\pi)$ there exist $\Delta\in SZ(\Lambda)$
such that
the strict transform  of $\Gamma$ by $\rho$ intersects
$\Delta$ and the value $e=\h{h}(\Delta)$ is the unique
one such that $I_z(\phi_{e},\Gamma)>I_{z}(\phi_w,\Gamma)$ for all $w\neq e$.
\end{theo}

A consequence of these results is Theorem \ref{sp/gen}:

\begin{theo}\label{sp/gen}
Let $\Phi _w$ be a  fiber of $\Lambda$. Then
the three following properties are equivalent:

\begin{enumerate}
\item $\Phi _{e}$ is a special fibre of $\Lambda$.
\item $I_z(\phi_e, C(\pi))>\displaystyle{\min}_{\phi \in \Lambda}I_z( \phi, C(\pi))$.
\item $\mu (\phi_{e})> \displaystyle{\min}_{\phi \in \Lambda}\mu (\phi )$.
\end{enumerate}

\end{theo}

In section 2 once we have set some preliminary results, we construct and study the minimal good resolution of $\Lambda$. In section 3, we prove Theorem \ref{theoprincipal} and
\ref{specialfibre} and in section 4 we show Theorem \ref{critloc}. To finish, in section 5,
we present some examples.

\section{Preliminary results and notations}

Let $(Z,z)$ be a normal surface singularity and let
$\rho: (X,E)\to (Z,z)$ be a good resolution of it. That means that $\rho$ is a
resolution of the singularity $(Z,z)$ such that the exceptional divisor
$E=\rho^{-1}(z)$ is a union of smooth projective
curves
$E=\bigcup_{\alpha\in G(\rho)} E_{\alpha}$ with normal crossings, in particular three
of them have empty intersection.

For $\alpha\in G(\rho)$ and for each
holomorphic function $f: (Z,z)\to (\C,0)$ let denote by $\nu_{\alpha}(f)$ the
vanishing order of
$\k f = f\circ \rho: X\to \C$ along  the irreducible exceptional curve
$E_{\alpha}$ ($\nu_{\alpha}$ is just the divisorial valuation defined by $E_{\alpha}$).
The divisor $(\k f)$
 defined by $\k f = f\circ \rho$ on $X$ could be written as
$$
(\k f) = (\w f) + \sum_{\alpha\in G(\rho)}\nu_{\alpha}(f) E_{\alpha}
$$
where, the local part $(\w f)$ is the strict transform of the germ $\{f=0\}$.
For each $\beta\in G(\rho)$ one has the known Mumford formula (see \cite{Mu}):
\begin{equation}
\label{mumford}
(\k f)\cdot E_{\beta} = (\w f)\cdot E_{\beta} + \sum_{\alpha} \nu_{\alpha}(f)
(E_{\alpha}\cdot E_{\beta}) = 0\; .
\end{equation}
(Here ``$\cdot$" stand for the intersection form on the smooth surface $X$).
Notice that the intersection matrix $(E_{\alpha}\cdot E_{\beta})$ is negative
definite and so $\{\nu_\alpha(f)\}$ is the unique solution of the linear system
defined by the equations (\ref{mumford}) above.

\subsection{Intersection multiplicity}

Let $C\subset (Z,z)$ be an irreducible germ of curve in $(Z,z)$ and let
$f\in \OO_{Z,z}$ be a function.
Let
$\varphi: (\C,0)\to (C,z)$ be a parametrization (uniformization) of $(C,z)$, then we
define the intersection multiplicity of $\{f=0\}\subset Z$ and $C$ at $z\in C$ as
$I_z(f,C)= \ord_{\tau}(f\circ \varphi(\tau))$ ($\tau$ is the parameter in $\C$).
Notice that the normalization $\k{\OO_{C,z}}$ of the local ring $\OO_{C,z}$ of the
germ $C$ at $z\in C$ is a discrete valuation ring, so $\k{\OO_{C,z}}\simeq \C\{t\}$
for a uniformizing parameter $t$ and the valuation $v_{C}$ is defined by the order
function on $t$, i.e. $v_{C}(g)= \ord_t(g(t))$ for $g\in \OO_{C,z}\subset \C\{t\}$.
One has also that  $I_{z}(f,C)=v_C(f)$.
The intersection multiplicity $I_{z}(f,C)$ could be also understood as
the degree
$\deg(f|C)$
of the composition map of $f|_{C}: C\setminus\{z\}\to \C^{*}$ and the map
from $\C^*$ into the unit circle $S^1$ which sends a non-zero complex number $t$,
onto $t/|t|$.
Obviously the above definition could be extended by linearity to define the
intersection multiplicity of a $f$ with a (local) divisor
$\sum_{i=1}^k n_i C_i$ as
$I_{z}(f,\sum n_i C_i) = \sum n_i I_z(f, C_i)$.

Let $\rho: (X,E)\to (Z,z)$ be a good resolution of the normal singularity $(Z,z)$ and
$E = \bigcup_{\alpha\in G(\rho)}E_{\alpha}$ be the exceptional
divisor.
Let $\w{C} := \k{\rho^{-1}(C\setminus \{z\})}$ be the strict transform of $C$ by
$\rho$. Then (see \cite{Mu})
$$
I_{z}(f,C) = (\k{f})\cdot \w{C} = (\w{f})\cdot \w C + \sum_{\alpha \in
G(\rho)}\nu_{\alpha}(f)
(E_{\alpha}\cdot \w{C})\; .
$$

Let us take now a good resolution $\rho$ such that $\w{C}$ is smooth and transversal
to $E$ at a smooth point $P$ and also with the condition $(\w{f})\cdot \w{C}=0$. This
resolution could be obtained by a finite number of point blowing ups starting on
(say) the minimal good resolution of $(Z,z)$. Let $\alpha(C)\in G(\rho)$ be the (unique)
component of $E$ such that $\w{C}\cap E_{\alpha(C)} = P$. Then one has
$I_{z}(f,C) = \nu_{\alpha(C)}(f) = I_{P}(f\circ \rho,\w{C})$. Here
$I_{P}(-,-)$ coincides with the usual local intersection multiplicity of two germs at the smooth
local surface $(X,P)$. Notice that $\w{C}$ is a curvetta at the point
$P\in E_{\alpha(C)}$,
$\w{C}$ is the normalization of $C$  and $\rho|_{\w{C}}: \w{C}\to C$ is a
uniformization of $C$.

Let $f,g$ be analytic functions on $(Z,z)$ and let
$\Lambda = \langle f, g \rangle = \{\phi_w = w_2 f - w_1 g \, | \, w=(w_1:w_2)\in
\C  \P^1\}$ be the pencil of analytic functions defined by $f$ and $g$.
As in the case of plane branches (see \cite{Ca}),
 one has the following easy and
useful result:

\begin{prop}\label{prop1}
Let $C\subset (Z,z)$ be an irreducible germ of curve. Then there exists a unique
$w_0\in \C  \P^1$ such that
$I_{z}(\phi_w,C)$ is constant for all $w\in \C  \P^1\setminus{w_0}$ and
$I_{z}(\phi_{w_0},C)> I_{z}(\phi_w,C)$.
\end{prop}

{\it Proof.} The statement is trivial taking into account that
the valuation defined by $C$, $\nu_{C}$, is the order of the series in
$\C\{t\}$.

\subsection{Resolution of pencils}

Let $\pi = (f,g): (Z,z)\to (\C^2,0)$ be finite complex analytic morphism germ, let
$\Lambda = \langle f, g \rangle = \{w_2 f - w_1 g \,|\, w=(w_1:w_2)\in \C  \P^1\}$ be
the pencil of analytic functions defined by $f$ and $g$ and
let $h = (f/g): V \to \C  \P^1$ be the meromorphic function defined by $f/g$ in a
suitable punctured neighbourhood of $z\in Z$.

A {\it good} resolution
of $(f,g)$ is a good resolution
$\rho: (Y,E)\to (Z,z)$
of $(Z,z)$ such that
the (reduced) divisor  $|(fg\circ \rho)^{-1}(0)|$ has normal
crossings. It means in particular that three irreducible components of
$|(fg\circ \rho)^{-1}(0)|$
doesn't meet at a same point. Starting on the
minimal good resolution of $(Z,z)$ one can produce a good
resolution of $(f,g)$ by a sequence of blowing-ups of points in the corresponding smooth
surface (essentially resolving the singularities of the reduced total transform of
the curve $\{fg=0\}$).
We also call it a good resolution of the corresponding curves $\Phi_{(0:1)}\cup
\Phi_{(1:0)}$.
Such a good resolution is  minimal if  and only if the contraction of any  rational
component of self-intersection -1 of the exceptional divisor
does not give  a good resolution anymore.

As  defined in the introduction, a {\it good} resolution
of the pencil
$\Lambda$ is a good resolution
$\rho: (X,E)\to (Z,z)$
of the singularity $(Z,z)$,
such that the lifting
$\h h = h\circ \rho$ is a morphism on $X$.

\bigskip
Let $\rho : (X,E)\to (Z,z)$ be a good resolution  of  $(Z,z)$ and $E_\alpha$ an irreducible component  of $E$.
The {\it Hironaka quotient} of $(f,g)$ on $E_\alpha$ is the following rational number:
$$q(E_\alpha):= \displaystyle\frac{\nu_{\alpha}(f)}{\nu_{\alpha}(g)}.$$

If  $q(E_\alpha)>1$ (resp.  $q(E_\alpha)<1$) then the component $E_\alpha$ belongs to the
zero divisor (resp. pole divisor) of $h\circ \rho$. Note that
if  $E_\alpha$ is a dicritical component  of $E$ then $q(E_\alpha)=1$.
Notice that there may exists irreducible components $E_\alpha$ of $E$ which are not dicritical and for
which $q(E_\alpha)=1$. Those are all components for which
the restriction of $h\circ \rho$ is constant on $E_\alpha$
and $E_\alpha$  does not belong to the zero divisor nor to the pole divisor.

\begin{prop}\label{propMGRP} There exists a (unique) minimal good resolution of $\Lambda$.
\end{prop}

{\it Proof.}  Let $\rho' :(Y', E')\to (Z,z)$  be the minimal good
resolution of $(f,g)$.
The indetermination points of $h\circ \rho'$ are the intersection points of
irreducible components  $E_{\alpha}$ and $E_{\beta}$ of the total transform
$| (fg\circ {\rho '})^{-1}(0)|$ for which one has  $q(E_\alpha)>1$ and
$q(E_\beta)<1$. Here one of the components, $E_{\alpha}$ or  $E_{\beta}$, is allowed to be
the strict transform $\w \xi$ of a branch $\xi$ of $\{f=0\}$ (in such a case we put
$q(\w{\xi})>1$)
or $\{g=0\}$ (respectively
$q(\w{\xi})<1$).
Let $P$ be such an indetermination point. Blowing-up at $P$ one creates a divisor
$E_\eta$ of genus $0$ and one has that
$\nu_{\eta}(f)=\nu_{\alpha}(f)+\nu_{\beta}(f)$ and
$\nu_{\eta}(g)=\nu_{\alpha}(g)+\nu_{\beta}(g)$.
(If $E_{\beta}$ is a branch $\xi$ of
$\{f=0\}$ of multiplicity $r$, we have $\nu_{\beta}(f)=r$ and $\nu_{\beta}(g)=0$. We use
similar conventions for the case in which $E_{\beta}$ is a branch of $\{g=0\}$.)
If
$q(E_\eta)=1$, then neither
$E_\alpha\cap E_\eta$ nor $E_\beta\cap E_\eta$ is an indetermination point and
moreover $E_{\eta}$ is a dicritical divisor. Else if
$ q(E_\eta)>1$ (resp. $ q(E_\eta)<1) $ then $E_\beta\cap E_\eta$ (resp.
$E_\alpha\cap E_\eta$)  is an indetermination point and  we iterate the process.
After a finite number of blow-ups there does not subsist indetermination points and
so we have constructed a good resolution
$\rho'':(Y'',E'')\to (Z,z)$
of $\Lambda$.

Now, to obtain a minimal good resolution of $\Lambda$, we
have to contract some rational component of self-intersection $-1$ of the exceptional
divisor (see theorem 5.9 of \cite{L1}). By the above construction the new components
(specially the last one which is dicritical and with self-intersection $-1$) can not
be contracted because in such a case we have an indetermination point. As a
consequence a minimal good resolution of $\Lambda$ is obtained from $\rho''$ by
iterated contraction of the rational component of self-intersection $-1$ of the
exceptional divisor which are not dicritical.
Uniqueness follows as in the case of the usual minimal resolution (see for example \cite{BPV} th. 6.2 p. 86).

\bigskip

Let consider  $\rho : (Y,E)\to (Z,z)$  the minimal good resolution of the pencil
$\Lambda$ and $\widehat h= h\circ \rho$. For $w\in \C\P^1$ let $\widehat{h}^{-1}(w) =
\w{\Phi_w}$ be the strict transform of the fibre $\Phi_w$.
For $D$
a dicritical  component of $E$,  we will denote by
deg$(\widehat h_{\mid D})$ the
degree of the restriction of $\widehat h$ to
$D$,
$\widehat{h}_{\mid D}: D\to \C\P^1$.

\begin{prop}\label{proptopgen}
Let $w$ be a generic value for the pencil $\Lambda$, then
\begin{itemize}
\item[a)] The resolution $\rho$ is a good resolution of $\phi_w$.
\item[b)] $\w{\Phi_w}$ intersects $E$ only at smooth points of ${\cal D}$.
\item[c)] If $D\in {\cal D}$, the number of intersection points of $\w{\Phi_w}$ and $D$ is equal to
deg$(\widehat{h}_{\mid D})$.
\end{itemize}
Moreover,
the minimal good resolution of $\Lambda$ is the minimal good resolution of any
pair  of generic elements of $\Lambda$.
\end{prop}

{\it Proof}. By definition of a generic value,  $\w{\Phi_w}$ meets the
exceptional divisor  $E$
only at smooth points of $\cal D$.
Let $D$ be an irreducible component of $\cal D$ and $P $ a point of $\w{\Phi_w}\cap D$. Then, as $P$ is
not a critical point for $\widehat h$,   $\w{\Phi_w}$ is smooth and transversal to $D$ at $P$.
This implies also that
$$
\mbox{deg}\ (\widehat h_{\mid D})= \displaystyle\sum_{P\in D}
I_{P}(\widetilde{\phi_w},D)
$$
So, one has deg$(\widehat{h}_{\mid D}) = \# (\widetilde{\Phi_w}\cap D)$.

Now, let
$w'$ be another generic value.
Notice that the strict transforms of
$\w{\Phi_w}$ and  $\w{\Phi_{w'}}$ intersect
in the same number of points each dicritical divisor $D$, so both fibres have the same number of
branches, just
$ \sum_{D\in {\cal D}} \mbox{deg}(\widehat{h}_{\mid D})$.
Moreover,
$\w{\Phi_w}$ and  $\w{\Phi_{w'}}$ do
not intersect $\cal D$ at the same points because $\widehat h$ is a morphism. As a
consequence the minimal good resolution of $\Lambda$ is a good resolution of any
pair of generic fibres.  It leaves to show that it is the minimal one.

By definition of the minimal good resolution of  $\Lambda$, the irreducible
components of the exceptional divisor of self-intersection $-1$  we have to
contract in  the minimal good resolution of the pencil  $\Lambda$,  to reach
the minimal good resolution of  $(\phi_w,\phi_{w'})$, lie  in the dicritical
components (see the proof of proposition \ref{propMGRP}; it is a consequence of the
construction of the minimal good resolution of $\Lambda$).
Contracting a dicritical component we obtain a map $\rho ''$ such that
the strict transforms of $\Phi_{w}$ and $\Phi_{w'}$ by $\rho''$
intersect an irreducible component of the
exceptional divisor at the same point and so
the strict transform by $\rho''$ of
$\{\phi_w\phi_{w'} =0\}$ has not
normal crossings. Consequently the minimal good resolution of $\Lambda$ is the
minimal good resolution of the pair $(\phi_w,\phi_{w'})$.

\subsection{Hironaka quotients}

In 2.2 we have defined the Hironaka quotient of $(f,g)$ on an
irreducible component  $E_\alpha$ of the exceptional divisor of a good resolution of $(Z,z)$.
In the same way we can define the Hironaka quotient of $(\phi_w,\phi_{w'})$  on
$E_\alpha$ for any pair $(\phi_w,\phi_{w'})$ of elements of $\Lambda =
\langle f, g \rangle$ as the rational number
$$q^w_{w'}(E_\alpha):=
\displaystyle\frac{\nu_{\alpha}(\phi_{w})}{\nu_{\alpha}(\phi_{w'})}.$$
In this way  $q(E_\alpha)= q^0_\infty (E_\alpha)$ (here $0=(0:1)\in \C  \P^1$,
$\infty=(1:0)\in \C  \P^1$)
but to simplify the notations we
will still write $q(E_\alpha)$ for the Hironaka quotient of $(f,g)$.

Notice that an irreducible component $E_\alpha$ of $E$ is dicritical if and only
if  $q^w_{w'}(E_\alpha)=1$ for any pair $(w,{w'})$ of elements of $\C \P^1$.

\begin{coro}\label{corQH}
The Hironaka quotient of any pair of generic elements of $\Lambda$ associated to
any irreducible component of the exceptional divisor of the minimal good resolution
of  $\Lambda$  is equal to one.
\end{coro}

{\it Proof}.
Let $w, w'\in \C\P^1$ be a pair of generic values of $\Lambda$ and $D\in {\cal D}$, then
$(\widetilde{\phi_w}) \cdot D = (\widetilde{\phi_{w'}}) \cdot D = \mbox{deg}(\widehat{h}_{\mid D})$ (see proposition \ref{proptopgen}). On
the other hand, if $E_{\beta}$ is a non-dicritical component of $E$ then one has
$(\widetilde{\phi_w}) \cdot E_{\beta} = (\widetilde{\phi_{w'}}) \cdot E_{\beta} = 0
$.
Now, the system of linear equations given by the formula (\ref{mumford}) for $\phi_w$ and
$\phi_{w'}$ is the same and so the solutions $\{\nu_{\alpha}(\phi_w)\}$ and
$\{\nu_{\alpha}(\phi_{w'})\}$ are the same. Thus, $\nu_\alpha(\phi_w)=\nu_{\alpha}(\phi_{w'})$
and $q^w_{w'}(E_\alpha)=1$
for any $\alpha\in G(\rho)$.

\bigskip

\begin{rema}
Let $E_{\alpha}$ be a non dicritical component of the exceptional divisor of the
minimal good resolution of the pencil $\Lambda$ and let $C$ be
a curvet in $E_{\alpha}$
(an irreducible smooth curve germ whose strict transform intersects $E_{\alpha}$
in a smooth point) such that $P= \w{C}\cap E_{\alpha}$ does not belong to
the strict transform of  any fibre $\Phi$  of $\Lambda$. One has $I_z(\phi, C) = \nu_{\alpha}(\phi)$  for any
$\phi\in \Lambda$ and by Proposition \ref{prop1}
there exists a unique $e\in \C\P^1$ such that
$\nu_{\alpha}(\phi_w)$ is constant for all $w\in\C\P^1\backslash \{e\}$ and
$\nu_{\alpha}(\phi_e) > \nu_{\alpha}(\phi_w)$. Moreover,
the above value $e\in\C\P^1$ must be a special value of $\Lambda$.
\end{rema}

\bigskip

Let $b: ( Z_{I},E_{I})\to (Z,z)$ be the normalized blow-up of the ideal $I=(f,g)$.
In \cite{B2} and \cite{BL} an element $\phi \in I $ is defined  to be
{\it general} if it is {\it superficial} and the strict transform of $\Phi = \{\phi=0\}$ by $b$
is smooth and
transverse to the exceptional divisor at smooth points. (See definition 2.1 of
\cite{B2}).
Proposition 2.2 of \cite{B2} allows to characterize general elements in terms of any
good resolution of $Z_I$, in particular one can use a good resolution
$\rho: (Y,E)\to (Z,z)$ of the pencil $\Lambda$. In this terms
one has that $\phi\in \Lambda$ is general if
$$
\nu_{\alpha}(\phi) = \nu_{\alpha}(I) = \min_{\phi\in I}\{\nu_\alpha(\phi)\} =
\min_{\phi\in \Lambda}\{\nu_\alpha(\phi)\}
$$
and moreover, the strict transform of  $\Phi$ by $\rho$ is smooth and
transversal to $E$.
By using the definition of the Milnor number of a germ of curve given in
\cite{Greuel_Buchweitz}, from  Theorem 1 and 2 of \cite{BL} one has that
$\phi\in \Lambda$ is general if and only if
$$
\mu (\phi) = \mu (I) := \min_{\phi\in I}\{\mu(\phi)\} =
\min_{\phi\in \Lambda} \{\mu(\phi)\}\; .
$$

Using proposition \ref{proptopgen} and the above results about Hironaka quotients we
have that $\Phi_{w}$ is a generic fibre
if and only
if $\phi_w$ is {\it general}. Moreover, one has also that
$\mu(\phi_{w})= \min_{\phi\in \Lambda}\{\mu(\phi)\}$ if and only if $\phi_{w}$ is
generic, so, $\mu(\phi_{w_0})> \min_{\phi\in \Lambda}\{\mu(\phi)\}$ if and only if
$w_{0}$ is a special value of $\Lambda$. This is the  equivalence of 1 and 3 in
Theorem 4.

\section{Topology of special fibres}

\subsection{Dual graph and topology}

Let $M:= Z \cap S^{2n-1}_\varepsilon$ where $S^{2n-1}_\varepsilon$ represents the boundary of the small
ball
of radius $\varepsilon$ of ${\C^n}$ centered at $z$. The manifold $M$ is called the {\it link} (see
\cite{Mu} and also \cite{Wahl})
of the
singularity $(Z,z)$.

Let $\phi_w$ be an element of $\Lambda$ and  $K_{\phi_w}:= \phi^{-1}_w(0)\cap M$. The {\it multilink
${\bf K}_{\phi_w}$ of $\phi_w$} is the oriented link $K_{\phi_w}$ weighted by the multiplicities of the
irreducible components of $\phi_w$. For $\varepsilon $ small enough, the topology of the multilink
${\bf K}_{\phi_w}$ in $M$ does not depend on the choice of~$\varepsilon$.

The fibres $\Phi _w$ and $\Phi_{w'}$ are said to be  {\it topologically equivalent}   if and only if
there exists a diffeomorphism of $M$ that send ${\bf K}_{\phi_w}$ on ${\bf K}_{\phi_{w'}}$ respecting
orientations and weights (see~\cite{LMW}).

\bigskip

Let $\rho _w: (X,E)\to (Z,z)$ be the minimal good resolution of $(Z,z)$ such that the divisor
$(\phi_w\circ \rho_w)$ has normal crossings. From Neumann (see \cite{N}), the topology of the multilink ${\bf
K}_{\phi_w}$ determines the minimal good resolution $\rho_w$,  where  the irreducible components of the
strict transform  of $\Phi _w$ by $\rho _w$ are weighted with their multiplicity and
taking into account the self-intersections and genus of the irreducible components of the exceptional divisor.
Conversely,  the Mumford formula (\cite{Mu})  and the fact that the intersection matrix $
(E_\alpha\cdot
E_\beta)_{\alpha, \beta \in G(\rho _w)}$ is negative definite (so invertible) imply that the  set
$\{\nu_\alpha (\phi _w), \alpha \in G(\rho_w)\}$ is uniquely defined and  so the divisor
$(\overline{\phi_w})$ (see section 2) is  uniquely determinate on $X$ from
the set $\{(\widetilde{\phi_w})\cdot E_{\alpha} \, | \, \alpha \in G(\rho_w)\}$.
As a consequence  the minimal
good resolution $\rho _w$ characterizes  the topology  of  the multilink ${\bf K}_{\phi_w}$.

Let $\rho: (X,E)\to (Z,z)$ be a good resolution of the normal surface singularity
$(Z,z)$, $E=\bigcup_{\alpha\in G(\rho)}E_{\alpha}$ its exceptional divisor.
It is useful to encode the information of the resolution $\rho$ by means of the so
called {\it dual graph} of $\rho$. The set of vertices of this graph is the set
$G(\rho)$, each vertex $\alpha$
is pondered by  $(\alpha,E_{\alpha}^2,g(E_{\alpha}))$ where  $E_{\alpha}^2$ is the self-intersection
of  $E_{\alpha}$,  and $g(E_{\alpha})$ its genus.
An intersection point between $E_\alpha$ and $E_\beta$ is
represented by an edge  linking the vertices $\alpha$ and $\beta$.

If we take $\rho$ as a good resolution of the local curve $C=\sum_{i=1}^{\ell}n_i
C_i$ (in particular if $C=\{\varphi=0\}$ for some function $\varphi$)
 one
add an arrow for each irreducible component $C_i$ of $C$ weighted by the
multiplicity $n_i$. In the case in which we deal with a good resolution
of pair of functions $(f,g)$, in the graph of $fg=0$ one mark with different colors
the arrows corresponding to branches of $\{f=0\}$ and those of $\{g=0\}$
(another possibility is to use different kinds of marks, say for example arrows for $f$
and stars for $g$). The sharp extremities of the arrows are considered as somekind of special vertices of the
graph. The notations ${\cal G}(\rho)$,
${\cal G}(\rho, \varphi)$ and  ${\cal G}(\rho, f, g)$ will be used for the dual graph
in each situation. Note that the case of a good resolution $\rho$ of the pencil
$\Lambda=\langle f, g\rangle$ is encoded by the dual graph
${\cal G}(\rho,\phi_w,\phi_{w'})$ for a pair of generic fibres.

\medskip
Following Neumann, one has:

\medskip

{\bf Statement:} The fibre $\Phi_{w}$ and $\Phi_{w'}$ are topologically equivalent if and only if the
graphs
${\cal G}(\rho_w)$ and  ${\cal G}(\rho_{w'})$ are the same.

\bigskip

Let  $\rho: (X,E)\to (Z,z)$ be a good resolution  of $(f,g)$ and let $E_\alpha$
be an irreducible component of $E$. We denote  ${\stackrel{\circ}{E}}_\alpha$ the
set of smooth points of
${E}_\alpha$ in the reduced total transform $|(fg\circ \rho)^{-1}(0)|$. An irreducible
component
${E}_\alpha$ (or its corresponding vertex $\alpha$) of
$E$ is a {\it rupture component}  if $\chi ({\stackrel{\circ}{E}}_\alpha) <0$, where
$\chi$ is the Euler characteristic. Note that
$\chi ({\stackrel{\circ}{E}}_\alpha)$ is equal
to $2-2g(E_{\alpha})-v(\alpha)$, where $v(\alpha)$ is the number of intersection
points of $E_{\alpha}$ with other components of the total transform of $fg=0$.
Thus, the rupture components are all  the rational ones with at least three different
edges or arrows and all the non-rational irreducible components.
We will say that
$\alpha$ is an {\it end} when $\chi(\stackrel{\circ}{E}_{\alpha})=1$.
Obviously $\alpha$ is an end if and only if $E_{\alpha}$ is rational and one has only one edge on
it.

\bigskip

The {\it neighbouring-set} of $E_\alpha$ in $X$ is the set constituted of
$E_\alpha$ union the irreducible components of the exceptional divisor and of the
strict transform of $\{fg=0\}$ that intersect $E_\alpha$. We denote it
$st (E_\alpha)$.

A {\it chain} of length $r$, $r\geq 3$, in $E$ is a connected part of $E$
constituted of a finite set of irreducible components $E_{\alpha _1}, \ldots ,
E_{\alpha _r}$ satisfying:
\begin{itemize}
\item
 $\chi ({\stackrel{\circ}{E_{\alpha_i}}})=0$, for $2\leq i \leq r-1$, and
\item
 $st (E_{\alpha _i})=\{  E_{\alpha _{i-1}}, E_{\alpha _i}, E_{\alpha _{i+1}}\}$
for $2\leq i\leq r-1$.
\end{itemize}

Notice that  the strict transform of $\{ fg=0\}$ does not intersect $\{ E_{\alpha_2},
\ldots , E_{\alpha _{r-1}}\}$.

A {\it cycle} of length $r$, $r\geq 3$, in $E$ is a chain such that
$st (E_{\alpha _r})=\{  E_{\alpha _{r-1}}, E_{\alpha _r}, E_{\alpha _1}\}$.
A {\it cycle} of length $2$ in $E$ is a connected part of $E$ constituted of  two
irreducible components $E_{\alpha _1},  E_{\alpha _2}$ such that  $\chi
({\stackrel{\circ}{E_{\alpha_2}}})=0$ and $st (E_{\alpha _2})=\{  E_{\alpha _{1}},
E_{\alpha _2}\}$.

The following result is a  direct generalization of proposition 1 and corollary 1 of
\cite{DM}.

\begin{prop}\label{lemmeDM}
Let $\rho: (X,E)\to (Z,z)$ be a good resolution of $(f,g)$. Let $E_\alpha$ be an
irreducible
component of the exceptional divisor such that the strict transform of $\{ fg=0\}$
does not intersect $E_\alpha$. Then there exists $E_\beta $ in $st (E_\alpha)$ such
that $q(E_\beta )>q(E_\alpha)$ if and only if there exists $E_\gamma$ in $st
(E_\alpha)$ such that $q(E_\gamma )<q(E_\alpha)$.

Moreover, if $\{ E_{\alpha _1}, \ldots , E_{\alpha _r}\}$, $r\geq 3$ is a chain,
then one of the following facts is true:
\begin{itemize}
\item[$\bullet$]
$q(E_{\alpha _i})<q(E_{\alpha_{i+1}})$ for
$1\leq i \leq r-1$.
\item[$\bullet$]
$q(E_{\alpha _i})>q(E_{\alpha_{i+1}})$ for
$1\leq i \leq r-1$.
\item[$\bullet$]
$q(E_{\alpha_i})$ is constant for $1\leq i \leq r$.
\end{itemize}

In particular,  if $E_{\alpha _r}$ is an end, then $q(E_{\alpha_i})$ is constant
for $1\leq i \leq r$ and
if $\{ E_{\alpha _1}, \ldots , E_{\alpha _r}\}$ is a cycle, then  $q(E_{\alpha_i})$
is constant for $1\leq i \leq r$.
\end{prop}

The proof almost repeats the proof of the refereed Proposition by using the
equations (\ref{mumford}) for $f$ and the divisor $E_{\alpha}$ as well as the same
equation for $g$.
As proposition \ref{lemmeDM} is a key result, we give back the proof for the first statement.

\bigskip
{\it Proof.}
By using equation (\ref{mumford})  for $f$ we have:
$$
0= (\k{f})\cdot E_{\alpha} = (\w{f})\cdot E_{\alpha} + \sum_{\gamma}\nu_{\gamma}(f)
(E_{\gamma}\cdot E_{\alpha}) =
\sum_{\eta\in st(E_{\alpha})} \nu_{\eta}(f) (E_{\eta}\cdot E_{\alpha})\; .
$$
The same equation is true for $g$ instead $f$ and thus one has:
\begin{equation}\label{eq_fg}
\begin{aligned}
\sum_{E_\eta\in st(E_\alpha), \eta\neq
\alpha}\nu_\eta(f)(E_\eta\cdot E_\alpha) & = (-E_{\alpha}^2) \; \nu_\alpha(f)
\\
\sum_{E_\eta\in st(E_\alpha), \eta\neq
\alpha}\nu_\eta(g)(E_\eta\cdot E_\alpha) & = (-E_{\alpha}^2)  \; \nu_\alpha(g)
\end{aligned}
\end{equation}

\medskip
Let suppose that $q(E_\eta) \geq  q (E_\alpha)$ for each $E_\eta \in
st(E_\alpha)$. This condition is equivalent to:
$$(E_\eta\cdot E_\alpha)
\nu_\eta(f)\nu_\alpha(g)\geq (E_\eta\cdot E_\alpha) \nu_\alpha (f)\nu_\eta (g)\, .
$$
As $q (E_\beta) > q (E_\alpha)$, we obtain:
$$
\nu_\alpha (g)\displaystyle\sum_{E_\eta\in st(E_\alpha), \eta\neq \alpha}(E_\eta
\cdot E_\alpha) \nu_\eta(f)>\nu_\alpha(f)\sum_{E_\eta\in st(E_\alpha), \eta\neq
\alpha}(E_\eta \cdot E_\alpha) \nu_\eta (g)\,.
$$
However, by using  the equations (\ref{eq_fg}), both sides of the above inequality are
equal to
$(-E_{\alpha}^2) \nu_{\alpha}(f)\nu_{\alpha}(g)$ and so
we reach a contradiction.

The others statements  of the proposition are direct consequences of  this result.

\subsection{Proof of Theorems 1 and 2}

Let $\rho : (Y,E)\to (Z,z)$ be the minimal good resolution of the pencil  $\Lambda$, $\widehat h= h\circ
\rho$.
If $w$
and $w'$ are  generic values for the pencil $\Lambda$, the Proposition \ref{proptopgen}, together with
the above Statement give

\begin{coro}
Let $w, w'\in \C\P^1$ be generic values of $\Lambda$.
Then, the  fibres $\Phi_{w}$ and $\Phi_{w'}$ are topologically equivalent.
\end{coro}

Thus, in order to finish the proof of Theorem
\ref{theoprincipal} it
only remains to show that a special  fibre $\Phi_{e}$ is
not topologically equivalent to a generic one.

Let $\Delta$ be an element of  $SZ(\Lambda)$ and $e = \widehat h (\Delta)$. We denote
$\Phi_e$
the fibre of $\Lambda$ associated to $e$ and by $\widetilde{\Phi_{e}}$ its strict transform by $\rho$.
The remaining part of Theorem \ref{theoprincipal} and Theorem \ref{specialfibre}
are direct consequences of the three following lemmas.

\begin{lemm} \label{lemmtype1}
If   $e $ is a  special value of $\Lambda$ associated  to a connected component $\Delta$  of
$\overline{E\backslash{\cal D}}$, then the strict transform of  $\Phi_e$ by $\rho$ intersects $\Delta$.
\end{lemm}

{\it Proof}.
Let us assume that $\widetilde{\Phi_{e}}\cap \Delta = \emptyset$.
Notice that if we change $\rho$ by
a good resolution of $\Lambda$ such that it is also a good resolution
of
$\Phi_e$ then the connected set $\Delta$ remains unchanged. So, we can keep the notations  we
use for $\rho$ for this new resolution.

Consider the Hironaka quotient with respect to $e$ and $w$ as a map
$q^e_w : E\to \Q$. Note
that
for any $E_\alpha $ in $ \Delta$,  we have $q^e_w(E_\alpha)>1$.
Let $E_\beta$ be an irreducible component  of  $\Delta$  such that $q^e_w (E_{\beta})\geq q^e_w (E_\alpha)$ for
each $E_\alpha$ in $\Delta$ and let $\Delta'$ be the maximal connected subset of $E$
such that $E_\beta \in \Delta '$ and $(q^e_{w})_{\mid \Delta'}$ is constant and equal to  $q^e_w (E_{\beta})$.
Notice that $E_{\beta}\subset \Delta'\subset \Delta$ because
$q^e_w(E_{\alpha})=1$ for any $E_{\alpha}$ such that
$E_{\alpha}\cap \Delta \neq \emptyset$ and $E_{\alpha}\not\subset \Delta$ (in fact such an $E_{\alpha}$
is a dicritical divisor).
Let now $E_{\gamma}\subset \Delta'$ and such that $st(E_{\gamma})\not\subset \Delta'$ and
$E_{\alpha}\in st(E_{\gamma})$, such that $E_{\alpha}\not\subset \Delta'$.
One has $q^e_w(E_{\beta}) > q^e_w(E_\alpha) >1$ if $E_{\alpha}\subset \Delta$ and
$q^e_w(E_{\beta}) > q^e_w(E_\alpha) = 1$ otherwise.
However, being $\Delta'\subset \Delta$, this contradicts Proposition \ref{lemmeDM} for the irreducible
component $E_{\gamma}$.

As a consequence
$\widetilde{\Phi_e}\cap \Delta\neq \emptyset$ and so $\Phi_e$ can not be topologically equivalent to
$\Phi_w$ for a generic value $w$.

\begin{lemm} \label{lemmtype2}
If   $e $ is a  special value of $\Lambda$ associated to a  smooth point $P$ of   $D$  in $ {\cal D}$
which is a critical point of  $\widehat h$, then the strict transform of $\Phi_e$ by $\rho$ intersects
$ D$ at $P$. Moreover it is not smooth and transversal to   $D$ at $P$.
\end{lemm}

{\it Proof}.
Blowing-up   at $P$ we create a divisor $E_\alpha$.
 As $P$ lies in the zero locus of $(\phi_e/\phi_{w})\circ \rho$, for any value $w\neq e$ we have
$q^{e}_w (E_\alpha)>1$. Moreover,  as $D$ is a dicritical component, $q^{e}_w (D)=1$.
Now, if we assume that $P\notin \widetilde{\Phi_e}$ then one can use Proposition
\ref{lemmeDM} for the new divisor $E_{\alpha}$ and we reach a contradiction.

Assume that $\widetilde{\Phi_e}$ is smooth and transversal to $D$ at the point $P$. In this case
we can choose local coordinates $\{u,v\}$ on $Y$
at $P$ in such a way that $\widetilde{\Phi_{w}} = \{v=0\}$ and $D =\{u=0\}$ on a neighbourhood $V$ of $P$.
So, the function $\phi_e\circ \rho$ is $u^{a}v$ on $V$ and, for a generic
value $w$, $\phi_w\circ \rho$ is $u^{b}\eta(u,v)$ for a unit $\eta$.
Note that $a = \nu_{D}(\phi_e) = \nu_{D}(\phi_w)=b$, being $D$ dicritical, and so
the expression of
$\widehat{h}$ at $P$ is $v\eta^{-1}(u,v)$.
Now, the restriction of $\widehat{h}$ to $D$ is given locally at $P$ as the map $v\mapsto v$.
Thus the point $P$ is not a critical (ramified) point of $\widehat{h}_{\mid D}: D\to \C\P^1$.

As a consequence $\widetilde{\Phi_e}$  is not smooth and transversal to $D$ at $P$, in particular
it can not be topologically equivalent to
$\Phi_w$ for a generic value $w$.

\begin{lemm}\label{lemmtype3}
If   $e $ is a  special value of $\Lambda$ associated to an intersection point $P$ between two
irreducible components  of $ {\cal D}$, then the strict transform of $\Phi _e$ by  $\rho$ intersects $ \cal D$ at $P$.
\end{lemm}

{\it Proof}.
Let $P =E_{\alpha_1}\cap E_{\alpha_2}$ such that  $E_{\alpha_1}$ and $ E_{\alpha_2}$  are dicritical
components.
Let us assume that $P\notin \widetilde{\Phi_{e}}$.
Blowing-up at $P$ we create a divisor  $E_\alpha$ satisfying
$\{ E_{\alpha_1}, E_\alpha,E_{\alpha_2}\} = st (E_\alpha)$.
As $q^{e}_w(E_{\alpha_1})=q^{e}_w(E_{\alpha_2})  =1$ and $q^{e}_w(E_{\alpha})>1$, we reach a
contradiction with proposition
\ref{lemmeDM}.

As a consequence, ${\Phi_e}$ is not resolved by $\rho$ and so could not be topologically
equivalent to a generic fibre $\Phi_w$.

\section{Behaviour of the critical locus}

Let $\pi=(f,g): (Z,z)\to (\C^2,0)$ be a finite complex analytic morphism. Following
Teissier (\cite{T}), the critical locus of $\pi$ is the analytic subspace
defined by the zeroth Fitting ideal $F_0(\Omega_{\pi})$ of the module $\Omega_{\pi}$
of relative differentials. The critical locus may have embedded components, however
we are only interested in the components of dimension one. So, we denote by $C(\pi)$
the divisorial part of the critical set with its non-reduced structure, i.e. each of
its components counted with its multiplicity, and we refer to $C(\pi)$ as the
critical locus of $\pi$. Note that out of the singular point $z\in Z$,  $C(\pi)$ is
defined by the vanishing of the jacobian determinant and also that $C(\pi)$ depends
on $\Lambda$ and not on the pair of functions of $\Lambda$ fixed to define the
corresponding finite morphism, so we denote it also by $C(\Lambda)$.
If we denote $\Gamma_i$, (resp. $n_i$) $i=1,\ldots, \ell$ the irreducible components
(branches) of
$C(\Lambda)$ (resp. their multiplicity) then $C(\Lambda)$ is the local divisor
$C(\Lambda)= \sum_{i=1}^{\ell} n_i \Gamma_i$.

\bigskip
Before proving  theorem \ref{critloc} and \ref{sp/gen}, let us first recall two
results from \cite{LMW} and \cite{Mi}.

\bigskip
Let  $(\phi_w,\phi_{w'})$ be  any pair of germs of the pencil $\Lambda$, let
$\rho ': (Y',E') \to (Z,z)$ be the minimal good resolution of $ (\phi_w,\phi_{w'})$, and
denote
$\Gamma (w,w'):= (\Gamma _k)_{k\in K}$  the set of irreducible components of
$C(\Lambda)$ which are not sent to a coordinate axe by $(\phi_w,\phi_{w'})$.
Let $Z_r$ be the set constituted of the union of the smooth points of $E'$ (smooth
points of $E'$ in $|(\phi_w\phi_{w'}\circ \rho)^{-1}(0)|$)  contained in an
irreducible
component of $E'$ with Hironaka quotients equal to $r$,  and  the intersection points
of two irreducible components of $E'$ of Hironaka quotient $r$.  The set $Z_r$ is
called the {\it $r$-zone} of $G(\rho ')$.
A connected component of $Z_r$ which contains at least one rupture vertex is called
a {\it $r$-rupture zone}.
Then from \cite{LMW} we have:

\bigskip
\noindent
{\bf Theorem A}. {\it The set
$\left\{\displaystyle\frac{I_z(\phi_{w}, \Gamma_k)}{I_z(\phi_{w'}, \Gamma_k)}, k\in
K\right\}$ is equal to the set of
Hironaka quotients associated to the rupture components of $G(\rho
',\phi_w,\phi_{w'})$.}

\bigskip
In \cite{Mi} a repartition in bunches of the branches of $\Gamma(w,w')$  is given as follows:

\bigskip
\noindent
{\bf Theorem B}.   {\it The intersection of the strict transform  of $\Gamma
(w,w')$ with a connected component of $Z_r$ is not empty if and only if it is a
$r$-rupture zone. Moreover if $\Gamma $ is an irreducible component of $\Gamma
(w,w')$ whose strict transform  intersects a  $r$-rupture zone then
$\displaystyle\frac{I_z(\phi_{w}, \Gamma)}{I_z(\phi_{w'}, \Gamma)}=r$.}

\bigskip

Next Lemma treats the case of irreducible components of the critical locus which are
also components of a fibre.

\begin{lemm}\label{multiplefibre}
Let $\xi$ be an irreducible component of a fibre $\Phi_e$, then $\xi$ is non
reduced if and only if $\xi$ is an irreducible component of $C(\Lambda )$.
\end{lemm}

{\it Proof}.
Let $\xi$ be an irreducible component of a fibre $\Phi_e$.
Let $w\in \C  \P^1$ be a generic value and
let $\rho': (Y',E')\to (Z,z)$ be
the minimal good resolution
of $(\phi_e ,\phi_w)$.
Let $\w{\xi}$ be the strict transform of $\xi$ by $\rho'$ and
let $P$ be the intersection point of $\widetilde \xi$ with the exceptional divisor
$E'$,
$P=\w \xi \cap E_{\alpha}= \w \xi \cap E'$.
We can choose a local system of coordinates $(u,v)$
in a neighbourhood $U\subset Y'$ of
$P=(0,0)$ such that
$u=0$ is an equation of $E_\alpha$, $v=0$ is a equation of $\w{\xi}$ and
the equation of $\w{\Phi_{\gamma}}$ at $P$ is
$u^a v^k$
where $a = \nu_{\alpha}(\phi_{e})$ and $k$ is the multiplicity of the
branch $\xi$ in $\Phi_{e}$.
On the other hand the equation of $\w{\Phi_w}$ at $P$ is $u^{b}\eta (u,v)$,
$b=\nu_{\alpha}(\phi_w)$ and $\eta(u,v)$ a unit. So, the expression of $\widehat{h}$
at  $P\in U$ is $u^{a-b}v^k (\eta(u,v))^{-1}$.

Let us first suppose that $\xi$ belongs to  $C(\Lambda )$. Let $Q $ be a point of
$\widetilde \xi \backslash \{ P\}$, say $Q$ has local coordinates
$(u_0, 0)$. The restriction of $\widehat{h}$ on a small disc $D(u_0,0)$ centered at
$Q$ in $u=u_0$ is $v^k \eta_0(u_0,v)$ with $\eta_0(u_0,v)$ a unit and $k>1$ because
$\xi$ lies in the ramification locus.
So, as $k$ is the multiplicity of $\xi$ in  ${\Phi_{e}}$,
$\xi$ is non reduced.

Conversely, if $\xi$ is an irreducible component of a fibre $\Phi_e$ which is
not reduced, the multiplicity   $k$  of $\xi$ in ${\Phi_{e}}$ satisfies  $k>1$.
Moreover the local equation of $\widehat{h}$
on any small   disc $D(t,0)$ centered at any point of local coordinates $(t,0)$ in
$U$ is $v^k \eta(t,v)$ with $\eta(t,v)$ a unit. As $k>1$, each point $(t,0)$ is a
ramification point and so $\w\xi$ lies in the ramification locus. Hence $\xi $  is
an irreducible component of $C(\Lambda )$.

\subsection{Proof  of theorem \ref{critloc} for singular points of ${\cal D}$ and critical points of
the restriction of $\widehat h $ to $\cal D$}

In the sequel $\rho: (Y,E)\to (Z,z)$ is the minimal good resolution of $\Lambda$ and
${\cal D}$ the dicritical locus of $E$.

\begin{prop}\label{criticallocus}
Let $P\in {\cal D}$ be such that $P\notin \overline{E\backslash {\cal D}}$. Then,
$P$ is a singular point of $\cal D$ or a critical point of $\h h _{|
{\cal D}}$ if and only if there exists an irreducible component $\Gamma$ of
$C(\Lambda )$ such that its strict transform intersects $\cal D$ at $P$. Moreover if
$\h h(P)=e$ then $I_z(\phi_e, \Gamma)>I_z(\phi _w,\Gamma)$ for all $w\neq e$.
\end{prop}

{\it Proof.}  Let us assume that
there exists an
irreducible component $\Gamma$ of $C(\Lambda)$ whose strict transform intersects
$\cal D$ at $P$.  Let $e= \h h (P)$ and denote by $D$ the irreducible component of
${\cal D}$ such that $P\in D$. If $\Gamma$ is a branch of $\Phi_e$ then it must be a
multiple irreducible component of it by the above Lemma and as a consequence the point $P$ is a
critical point of $\h{h}|{\cal D}$.

So, let us consider the case in which
$\Gamma$ is not a branch of $\Phi_{e}$ and assume that $P$ is not a singular point of
${\cal D}$, i.e. $P$ is a smooth point of ${\cal D}$ in the exceptional divisor
$E$.

If the strict transform $\w{\Phi}_e$ of $\Phi_e$ at $P$ has normal crossings with
${\cal D}$, then there exists an irreducible branch $\xi$ of $\Phi_e$ such that
its strict transform $\w{\xi}$ coincides with $(\w{\Phi}_e)_{P}$, i.e.
$\w{\xi}$
is smooth, transversal to $D$ and $\xi$ is not a multiple branch of $\Phi_e$ by Lemma
\ref{multiplefibre}.
By Theorem B there exists a $r$-rupture zone $R$ in the minimal good resolution
of
$(\phi_e,\phi_w)$ (here $w$ is assumed to be a generic value) such that
the strict transform of $\Gamma$ intersects $R$ and moreover
$I_z(\phi_e,\Gamma)/I_z(\phi_w,\Gamma) = r$ with $r>1$ because $P\in \w{\Gamma}\cap \w{\Phi_e}$.
Taking into account that $\w{\Phi_e}$
is smooth and transversal to the dicritical divisor $D$, then
one has that
$P = \w{\Gamma}\cap E \subset D\subset R$ and so, by Theorem A,
$$
\frac{I_z(\phi_e,\Gamma)}{I_z(\phi_w,\Gamma)} = q^{e}_{w}(D)
= \frac{\nu_{D}(\phi_e)}{\nu_{D}(\phi_w)}\;.
$$
However this is impossible because the last quotient is equal to 1, being $D$
dicritical.
Thus, as a consequence, $(\w{\Phi_e})_P$ must be singular or tangent to $D$. In both
cases $P$ is a critical point of $\h{h}|_{\cal D}$ (i.e. $\phi_e$ is a special
function of $\Lambda$).

\bigskip
Conversely, let $P$ be a singular point of $\cal D$ or a smooth point of $\cal D$
which is a critical point of  $\h h_{| {\cal D}}$ and let $e=\widehat h(P)$, then
from  Theorem 2, $\Phi _e$ is a special fibre of $\Lambda$. If the
irreducible component of  $\widetilde\Phi _{e}$ that intersects $\cal D$ at $P$ is
non reduced then from lemma \ref{multiplefibre} we have finished. Thus, we assume
that $\widetilde \Phi_e$ is reduced at $P$.

First, note that $\w{\Phi}_e$ has not normal crossings with $E$ at $P$.
 Because  if $P$ is smooth on ${\cal D}$ then $\w{\Phi}_e$ is either singular or tangent to
$\cal D$, and in the other case, it means if $P$ is a  singular point of $\cal D$, then  there are at least three components of the total transform
intersecting at $P$.

Let $w,w'\in \C  \P^1$ be  generic values and
let $\rho': (Y',E')\to (Z,z)$ be the minimal good resolution of
$\phi_{w}\phi_{w'}\phi_{e}$. Note that $\rho' = \rho \circ \sigma$, where
$\sigma$ is a sequence of point blowing-ups on $Y$, each of them
produces some new irreducible rational exceptional components. In particular
$\Delta = \sigma^{-1}(P)\subset E'$ is a connected exceptional part and must contain
a rupture component $E_{\alpha}\subset E'$.
Notice that no component of $\Delta$ is contracted in the minimal good resolution
$\rho'': (Y'',E'')\to (Z,z)$ of the pair $(\phi_e, \phi_w)$; i.e
$\Delta\subset E''$ and in particular $E_{\alpha}\subset E''$ is also a rupture
component in $E''$.
Let $R$ be the corresponding rupture zone in $E''$ which contains $E_{\alpha}$.
Note that for each $E_{\beta}\subset R\subset \Delta$ one has
$r= q^e_{w}(E_{\beta}) = \nu_{\beta}(\phi_e)/\nu_{\beta}(\phi_w)>1$.

Now, Theorem A implies that there exists a branch $\Gamma$ of
$C(\Lambda)$ such that its strict transform by $\rho''$ intersects $\Delta$ and also
$$
\frac{I_{z}(\phi_e,\Gamma)}{I_{z}(\phi_w,\Gamma)} =
\frac{\nu_{\alpha}(\phi_e)}{\nu_{\alpha}(\phi_w)} = r >1.
$$
Taking into account that $R\subset \Delta$ and
$\sigma(\Delta)=P$,
 one has that
the strict transform of $\Gamma$ by $\rho$ intersects $E$ at the point $P$ and
moreover
$I_{z}(\phi_e,\Gamma) > I_{z}(\phi_w,\Gamma)$. Note that the above inequality is true
for any irreducible component $\Gamma$ of $C(\Lambda)$ such that
its strict transform by $\rho$ intersects $\cal D$ at P. Thus,
the special fibre $\phi_e$ is the unique fibre with the condition
$I_{z}(\phi_e,\Gamma) > \min_{w}I_{z}(\phi_w,\Gamma)$.

\bigskip
\begin{rema}
Notice that if  $P$ is a smooth point of $\cal D$ which is a critical point of
$\widehat h _{| {\cal D}}$ or if  $P$ is a singular point of $\cal D$, then for any
fibre $\Phi _a$ and $\Phi _{a'}$ different from $\Phi_{e}$, we have
$q^{a}_{a'}(E_{\alpha}) = 1$ and then
$I_{z}(\phi_{a},\Gamma)=I_{z}(\phi_{a'}, \Gamma)$.
\end{rema}

\subsection{Proof of theorem \ref{critloc} for the connected components of $\overline{E\backslash
{\cal D}}$}

\bigskip

 Let  us remind that $\rho : (Y,E)\to (Z,z)$ is the minimal good resolution of
$\Lambda$ and  ${\cal D}$ the dicritical locus of $ E$.
Let $\Delta $ be a connected component  of $\overline{E\backslash {\cal D}}$ such
that $(h\circ \rho )(\Delta ) = e$.
Let $w, w'$ be generic values of $\Lambda$ and let us denote
$\rho' : (Y',E')\to (Z,z)$ the minimal good resolution of
$\phi_w \phi_{w'} \phi_e$.
Let us denote by $\tau : (Y',E')\to (Y,E)$
the composition of point blowing-ups which produces $Y'$ from $(Y,E)$
\begin{equation*}
(Y',E')  \overset{\tau}{\to}  (Y,E) \overset{\rho}{\to} (Z,z)
\end{equation*}

Let $\Delta'$ by the pull-back of $\Delta$ by $\tau$. Note that $\Delta'$ is a connected component of
$\overline{E'\backslash {\cal D'}}$ because the dicritical locus ${\cal D'}$ on $E'$ is just the
strict transform of ${\cal D}$ by $\tau$.
We will distinguish two cases, depending on the
existence of a rupture component $E'_{\alpha}$ in $\Delta'$ (with respect  to $\phi_w$ and
$\phi_{e}$).

\bigskip

{\bf Case 1)} There exist a rupture component $E'_{\alpha}$ in $\Delta'$.

For each component $E_{\beta}\subset \Delta'$ one has
$q_{w'}^{w}(E_{\beta})=1$ and $q^e_{w}(E_{\beta})>1$. Let $R$ be the rupture zone
of $E'$ such that $E_\alpha\subset R$. Then $R\subset \Delta'$ because
$q^e_{w}$ is constant and $> 1$ on $R$ and moreover
$q^{e}_{w}(D)=1$ for any dicritical divisor, in particular for $D$ dicritical such
that $D\cap \Delta'\neq \emptyset$.

Now, from Theorem B, there exist a branch $\Gamma$ of the critical locus $C(\Lambda)$ such that its
strict transform by $\rho'$, $\widetilde{\Gamma}$, intersect $R$. As consequence the strict transform of
$\Gamma$ by $\rho$, $\tau(\widetilde{\Gamma})$ intersects $\Delta$.
Again Theorem B implies that
$q^{e}_{w}(E_\alpha) = I_{z}(\phi_{e}, \Gamma)/I_{z}(\phi_{w}, \Gamma)$ and so
the special value $e$ is the unique one such that
$I_{z}(\phi_{e}, \Gamma) > I_{z}(\phi_{w'}, \Gamma)$ for any generic value $w'$.

\bigskip

{\bf Case 2)} There are no rupture components in $\Delta'$.

In this case $\Delta' = \{E_{\alpha_1}, \ldots, E_{\alpha_r}\}$ in such a way that there exists
a dicritical component $D \in { \cal D'}$ such that
$\{D = E_{\alpha_0}, E_{\alpha_1}, \ldots, E_{\alpha_r}\}$ is a chain and $\chi(E_{\alpha_r})\ge 0$.
Now, note that the strict transform of $\Phi_{e}$ intersects $\Delta'$ (see Theorem \ref{specialfibre}),
so the only way to avoid the existence of a rupture component with respect to $\phi_w \phi_e$ is that
$E_{\alpha_r}$ is an end (i.e it is rational and is connected only with the previous one
$E_{\alpha_{r-1}}$) and such that $\widetilde{\Phi_{e}}$, the strict transform of $\Phi_{e}$
by $\rho'$, intersect $E_{\alpha_r}$. Moreover, $\widetilde{\Phi_{e}}$ with its reduced structure
is smooth and transversal to $E_{\alpha_r}$.
It means that the minimal good resolution of $\Lambda$ is a resolution of
the reduced irreducible component $\xi_{e}$ of $\Phi_{e}$ whose strict
transform meets $\Delta$  at $E_{\alpha_r}$. Actually, otherwise to resolve
$\xi_{e}$,  we have to blow-up at $\xi_{e}\cap E_{\alpha_r}$ and this process
produces a rupture component.

\begin{figure}[h]
$$\unitlength=1.00mm
\hspace*{-30mm}\begin{picture}(100.00,20.00)(0,10)
\thinlines
\put(10,15){\ldots\ldots}
\put(20,15){\vdots}
\put(20,23){\vdots}
\put(20,5){\vdots}
\put(20,10){\vdots}
\put(20.5,15){\circle{3}}
\put(20.5,15){\circle*{1.5}}
\put(20,15){\line(1,0){20}}
\put(22,10){$E_{\alpha_0}$}
\put(40,15){\circle*{1.5}}
\put(40,10){$E_{\alpha_1}$}
\put(40,15){\line(1,0){10}}
\put(50,15){\ldots\ldots}
\put(70,15){\line(1,0){40}}
\put(70,15){\circle*{1.5}}
\put(90,15){\circle*{1.5}}
\put(92,10){$E_{\alpha_{r-1}}$}
\put(110,15){\circle*{1.5}}
\put(112,10){$E_{\alpha_r}$}
\put(110,15){\vector(1,1){7}}
\end{picture}
$$
\caption{
Graph in Case 2
} \label{fig0}
\end{figure}
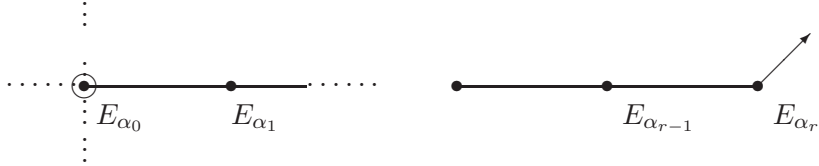

\begin{lemm}\label{lemtec}
Let $\seq{v}0r$, $\seq{e}1r$ be sequences of integers
such that
$v_{i-1} = e_i v_i -v_{i+1}$ for $i=1,\ldots, r-1$. Let
$\seq{q}0{r-1}\in \Z$ defined recursively as
$q_0=1$, $q_1=e_1$ and, for $i\ge 2$,  $q_{i}=e_iq_{i-1}-q_{i-2}$.
Then, for $i\ge 1$ one has $\gcd(q_i, q_{i-1})=1$ and
$v_0=q_{i} v_{i}-q_{i-1}v_{i+1}$.
\end{lemm}

{\it Proof.}
Obviously $\gcd(q_0,q_1)=1$ and from the definition of $q_{i}$,
if $\gcd(q_{i-1},q_{i-2})=1$ then $\gcd(q_{i-1},q_i)=1$. The
equality $v_0=q_{i}v_i-q_{i-1}v_{i+1}$ is obvious for $i=1$ and,
by induction, using the equality
$v_{i-1} = e_i v_i -v_{i+1}$ in the inductive hypothesis
$v_0=q_{i-1}v_{i-1}-q_{i-2}v_{i}$ one has
$$
v_0= q_{i-1}v_{i-1}-q_{i-2}v_{i} = q_{i-1}(e_i v_i -v_{i+1})
-q_{i-2}v_{i} = q_{i}v_{i}-q_{i-1}v_{i+1}\; .
$$

Now, the proof of the case 2 is a consequence of the next:

\begin{prop}
The irreducible curve  $\xi_{e}$ is a branch of $\Phi_{e}$ with multiplicity
bigger than 1.
As a consequence $\xi_e$ is also a branch of $C(\Lambda)$ and
so $C(\Lambda)$ intersect  $\Delta$.
\end{prop}

{\it Proof.} Recall that
$w$ is a generic element
of $\Lambda$.  For the sake of simplicity let denote
$v_i=\nu_{\alpha_i}(\phi_w)$ and $e_i = - E_{\alpha_i}^2$ for $i=0,\ldots,r$.
Then, by using the formula
\begin{equation}\label{intersect}
\left(
(\widetilde{\phi_{w}}) + \sum_{\alpha\in G(\rho')} \nu_{\alpha}(\phi_w) E_{\alpha}
\right) \cdot E_{\alpha_i} = 0
\end{equation}
for $i=1,\ldots, r$ one has that
\begin{equation}\label{dead}
\begin{aligned}
v_0 &= e_1 v_1 - v_2\\
v_1 &= e_2 v_2 - v_3 \\
 & \cdots \\
v_{r-2} &= e_{r-1} v_{r-1} - v_r\\
v_{r-1} &= e_r v_r
\end{aligned}
\end{equation}

By Lemma \ref{lemtec} one has $v_0 = q_rv_r$. Moreover, taking into account that
$e_i = - E_{\alpha_i}^2 \geq 2$ one can easily prove that
$q_r > q_{r-1} > \cdots > q_1 > q_0 = 1$.

Let us consider now the special fibre $\Phi_{e}$ and
let us denote $v'_i=\nu_{E_{\alpha_i}}(\phi_{e})$ for $i=0,\ldots, r$.
The equations (\ref{intersect}) applied for $\phi_e$ instead of $\phi_w$ gives a sequence of
equalities
$v'_{i-1} = e_i v'_i -v'_{i+1}$, for $i=1,\ldots, r-1$  (like in (\ref{dead}) above
with $v'_i$ instead $v_i$) together with the last one:
$$
v'_{r-1} = e_r v'_{r} - (\widetilde{\phi_{e}})\cdot E_{\sigma}
=e_r v'_{r} - k \; .
$$
Lemma \ref{lemtec} implies that $v'_0 = q_r v'_r - q_{r-1} k$.
Being $E_{\alpha_0}=D$ a dicritical divisor one has that
$v'_0  = \nu_{\alpha_0}(\phi_{e}) =
\nu_{\alpha_0}(\phi_{w}) = v_0$, i.e.
$$
q_r v_r = q_r v'_r -q_{r-1} k \; .
$$
By Lemma \ref{lemtec} again,
$\gcd(q_r,q_{r-1})=1$ and so $q_r$ divides $k$. In particular $k =
(\widetilde{\phi_{e}})\cdot E_{\sigma}>1$
and the irreducible germ $\xi_e$ appears repeated $k$ times
in
$\Phi_{e}$.

\subsection{Special fibres and critical locus}

Let $C(\Lambda)= \sum_{i=1}^{\ell}n_i \Gamma_i$ be the decomposition of
the critical locus in irreducible components.
For each $i\in \{1,\ldots, \ell\}$
the intersection multiplicity
$I_z(\phi, \Gamma_i)$ is constant but for exactly the unique
special value $\varepsilon(\Gamma_i) (=\varepsilon(i))$
such that
$I_z(\phi_{\varepsilon(i)}, \Gamma_i) >
I_z(\phi, \Gamma_i)$, for $\phi\neq \phi_{\varepsilon(i)}$.
So, as in  \cite{DM}, one has a
surjective map $\varepsilon : {\cal B}(C(\Lambda))\to Sp(\Lambda)$ from the set of
branches of the
critical locus to the set of special fibres of $\Lambda$.

If $w\in \C\P^1$ is a generic value one has that
$$
I_{z}(\phi_w, C(\Lambda)) = \sum_{i=1}^{\ell}n_i I_z(\phi_w,\Gamma_i)  =
\min \{ I_z(\phi, C(\Lambda)) \, , \, \phi\in \Lambda\}
$$
and, on the other hand, for a special value
$e\in \C\P^1$ one has
$$
I_{z}(\phi_e, C(\Lambda)) = \sum_{i=1}^{\ell}n_i I_z(\phi_e,\Gamma_i)  >
\sum_{i=1}^{\ell}n_i I_z(\phi_w,\Gamma_i)  =
\min \{ I_z(\phi, C(\Lambda)) \, , \, \phi\in \Lambda\}\; .
$$
Thus, as a consequence one has the following

\begin{coro}
$\Phi_e$ is a special fibre of $\Lambda$  if and only if
$$
I_z (\phi _e, C(\Lambda) >  \min  \left\{
I_z(\phi , C(\Lambda) , \phi \in \Lambda\right\}\; .
$$
\end{coro}

\begin{rema}
As in \cite{DM} the map
$\varepsilon: {\cal B}(C(\Lambda))\to Sp(\Lambda)$, defined above,  could be factorized
through the set of special zones $SZ(\Lambda)$
as $\varepsilon = \xi\circ \psi$:
$$
{\cal B}(C(\Lambda))\overset{\psi}{\to} SZ(\Lambda) \overset{\xi}{\to} Sp(\Lambda)
$$
The map $\psi$ associates to the branch $\Gamma$ the special zone
$\Delta$ such that the strict transform of $\Gamma$ in the minimal good resolution
intersects $\Delta$. In the same way the map $\xi$ sends $\Delta\in SZ(\Lambda)$
to $\widehat{h}(\Delta)$.

By means of a good resolution of all the fibres of $\Lambda$
$\rho' : (Y',E')\to (Z,z)$ (i.e. a good resolution of the product of all the
special fibres and a pair of generic ones) and the determination of all the rupture
zones in $E'$ with respect to the pairs $(\phi_e,\phi_w)$, being $e$ special and $w$
generic, one can determine a finer decomposition in bunches of the branches of the
critical locus $C(\Lambda)$.
\end{rema}

\section{Examples}

As seen in section 3.1,   to the minimal good resolution $\rho$ of the pencil $\Lambda$, one can
associate its intersection graph ${\cal G}(\rho)$.
The following examples illustrate  theorems 1,  2 and 3  in terms of
intersection graph. To construct ${\cal G}(\rho)$, we
follow the method of Laufer described in \cite{L2}, \cite{LW2} and also \cite{LMW}. It consists in
first establishing
the graph of  the minimal resolution of the discriminant curve,  which is the  image by $\pi$ of the
critical locus $C(\pi)$ of $\pi$.
Then we deduce the graph of the minimal good resolution of $(Z,z)$ and then the one of $G(\rho)$,
using  in particular   proposition 3.6.1 and 3.7.1 of \cite{LW2}.
As in the Figure~\ref{fig0} of Section 4 we use a
different kind of mark for the vertices representing dicritical divisors.

 \subsection{Example 1}

 Let $(Z,z)$ be defined by
 $z^3 =h(x,y)$
 with $ h(x,y)= (y+x^2)(y-x^2)(y+2x^2)(x+y^2)(x-y^2)(x+2y^2)$ and let $\pi$ be the projection on the
$(x,y)$-plane. Such a way $(u,v)=(x,y)$ and $f=u\circ \pi = x$ and $g=v\circ \pi=y$.

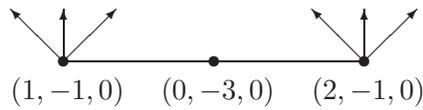
\begin{figure}[h]
\unitlength=1.00mm
$$
\hspace*{-30mm}\begin{picture}(100.00,20.00)(0,10)
\thinlines

\put(50,15){\circle*{1.5}}
\put(43,10){$(1,-1,0)$}
\put(50,15){\line(1,0){40}}
\put(70,15){\circle*{1.5}}
\put(63,10){$(0,-3,0)$}
\put(90,15){\circle*{1.5}}
\put(83,10){$(2,-1,0)$}
\put(50,15){\vector(-1,1){7}}
\put(50,15){\vector(0,1){7}}
\put(50,15){\vector(1,1){7}}
\put(90,15){\vector(1,1){7}}
\put(90,15){\vector(-1,1){7}}
\put(90,15){\vector(0,1){7}}

\end{picture}
$$
\caption{
Graph of  the discriminant of $\pi$.
} \label{fig1}
\end{figure}

The discriminant curve  of $\pi$ is  the curve $h(u,v)=0$.
The dual graph of its minimal embedded resolution is represented in the Figure \ref{fig1}.

\begin{figure}[h]
\unitlength=1.00mm
$$
\hspace*{-30mm}\begin{picture}(100.00,20.00)(0,10)
\thinlines

\put(50,15){\circle*{1.5}}
\put(38,10){$(1, -3, 1)$}
\put(50,15){\line(1,0){40}}
\put(70,5){\circle*{1.5}}
\put(70,15){\circle*{1.5}}
\put(70,25){\circle*{1.5}}
\put(63,17){$(0^2,-3, 0)$}
\put(63,0){$(0^1,-3, 0)$}
\put(63,28){$(0^3,-3, 0)$}
\put(90,15){\circle*{1.5}}
\put(88,10){$(2,-3,1)$}
\put(50,15){\line(2,1){20}}
\put(50,15){\line(2,-1){20}}
\put(90,15){\line(-2,1){20}}
\put(90,15){\line(-2,-1){20}}

\end{picture}
$$
\bigskip
\caption{The graph of  the minimal good resolution of $(Z,z)$.
} \label{fig2}
\end{figure}
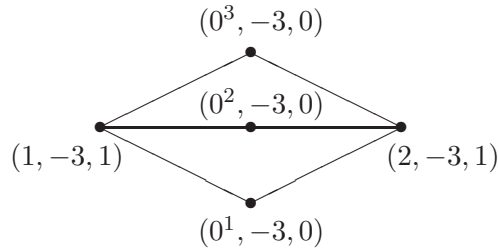

From proposition 3.6.1 of \cite{LW2} we deduce the graph of the minimal good resolution of~$(Z,z)$ (see
Figure \ref{fig2}).

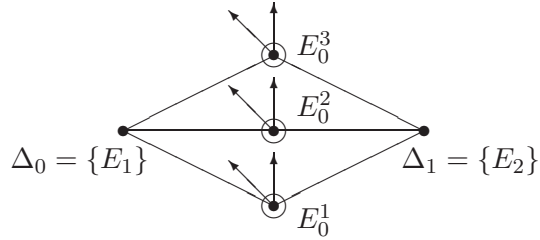
\begin{figure}[h]
\unitlength=1.00mm
$$
\hspace*{-30mm}\begin{picture}(100.00,20.00)(0,10)
\thinlines

\put(50,15){\circle*{1.5}}
\put(50,15){\line(1,0){40}}
\put(70,5){\circle*{1.5}}
\put(70,15){\circle*{1.5}}
\put(70,25){\circle*{1.5}}
\put(70,5){\circle{3}}
\put(70,15){\circle{3}}
\put(70,25){\circle{3}}
\put(90,15){\circle*{1.5}}
\put(50,15){\line(2,1){20}}
\put(50,15){\line(2,-1){20}}
\put(90,15){\line(-2,1){20}}
\put(90,15){\line(-2,-1){20}}
\put(70,5){\vector(-1,1){6}}
\put(70,15){\vector(-1,1){6}}
\put(70,25){\vector(-1,1){6}}
\put(70,5){\vector(0,1){7}}
\put(70,15){\vector(0,1){7}}
\put(70,25){\vector(0,1){7}}
\put(73,2){$E_{0}^1$}
\put(73,17){$E_{0}^2$}
\put(73,25){$E_{0}^3$}
\put(35,10){$\Delta_0=\{ E_1\}$}
\put(87,10){$\Delta_1=\{ E_2\}$}

\end{picture}
$$
\bigskip
\caption{The graph of  the minimal good resolution of $\Lambda$.
} \label{fig3}
\end{figure}

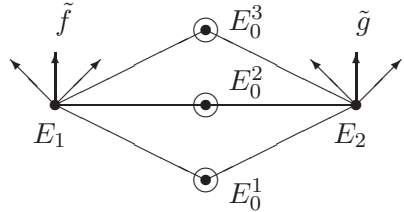
\begin{figure}[h]
\unitlength=1.00mm
$$
\hspace*{-30mm}\begin{picture}(100.00,20.00)(0,10)
\thinlines
\put(50,15){\circle*{1.5}}
\put(50,15){\line(1,0){40}}
\put(70,5){\circle*{1.5}}
\put(70,15){\circle*{1.5}}
\put(70,25){\circle*{1.5}}
\put(70,5){\circle{3}}
\put(70,15){\circle{3}}
\put(70,25){\circle{3}}
\put(90,15){\circle*{1.5}}
\put(50,15){\line(2,1){20}}
\put(50,15){\line(2,-1){20}}
\put(90,15){\line(-2,1){20}}
\put(90,15){\line(-2,-1){20}}
\put(50,15){\vector(0,1){7}}
\put(50,15){\vector(-1,1){6}}
\put(50,15){\vector(1,1){6}}
\put(50,25){$\tilde f$}
\put(90,15){\vector(0,1){7}}
\put(90,15){\vector(-1,1){6}}
\put(90,15){\vector(1,1){6}}
\put(90,25){$\tilde g$}
\put(73,2){$E_{0}^1$}
\put(73,17){$E_{0}^2$}
\put(73,25){$E_{0}^3$}
\put(47,10){$E_1$}
\put(87,10){$E_2$}

\end{picture}
$$
\bigskip
\caption{Minimal good resolution  of $(f,g)$.
} \label{fig4}
\end{figure}

As the minimal embedded resolution of the discriminant curve $h(u,v)=0$ of $\pi$ is also the minimal
good resolution of the product $uv(\lambda u+\mu v)h(u,v)=0$, for $(\lambda:\mu)\in \CP^1$,  from
propositions 3.6.1 and 3.7.1 of \cite{LW2}
we can deduce the graph of the minimal good resolution of $\Lambda$ (Figure \ref{fig3}), the one of
$(f,g)$ and as a consequence the one of  the minimal good resolution  of $(\phi_w\phi_{w'}fg)^{-1}(0)$
where $w$ and $w'$ are generic values of $\Lambda$ (Figure \ref{fig4}). Notice that the minimal good
resolution of $\Lambda$ is also the minimal good resolution of $(f,g)$.

The dicritical components of $E$ are $E_0^1, E_0^2, E_0^3$.  We have
$SZ(\Lambda)=\{\Delta_1,\Delta_2\}$ with  $\Delta_1= \{E_1\}$ and  $\Delta_2=\{ E_2\}$. The map
$(f/g)\circ \rho$ has no critical point  on $\cal D$ and $\cal D$ has no singular point neither. The
special fibre associated to $\Delta _1$ is  $\{ f=0\}$ and the one associated to $\Delta_2$ is $\{
g=0\}$.  We conclude that $\Lambda$ admits two special elements $f$ and $g$; the special value
associated to $\Delta_1$ is  $(0:1)$ and the one associated to $\Delta_2$ is $(1:0)$. The
Hironaka quotients are $q(E_1)=2$ and $q(E_2)=1/2$.

Moreover, using the minimal resolution of the discriminant curve (see
Figure~\ref{fig1}), we deduce  that,
for each $\Delta _i$,  there exists  three  irreducible components of  the reduced critical locus of $\pi$ whose strict transform intersects $\Delta_i$.

 \subsection{Example 2}

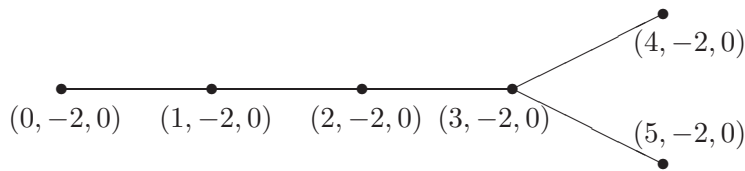
\begin{figure}[h]
$$\unitlength=1.00mm
\hspace*{-80mm}
\begin{picture}(100.00,20.00)(0,10)
\thinlines

\put(50,15){\circle*{1.5}}
\put(43,10){$(0,-2,0)$}
\put(50,15){\line(1,0){60}}
\put(70,15){\circle*{1.5}}
\put(63,10){$(1,-2,0)$}
\put(90,15){\circle*{1.5}}
\put(83,10){$(2,-2,0)$}
\put(110,15){\circle*{1.5}}
\put(130,25){\circle*{1.5}}
\put(100,10){$(3,-2,0)$}
\put(130,5){\circle*{1.5}}
\put(126,20){$(4,-2,0)$}
\put(126,8){$(5,-2,0)$}
\put(110,15){\line(2,1){20}}
\put(110,15){\line(2,-1){20}}

\end{picture}
$$
\caption{The graph of  the minimal good resolution of $D_6$.
} \label{fig2.2}
\end{figure}

 Let $(Z,z)$ be the singularity $D_6$ defined by the  equation
$z^2 =y(x^2+y^4)$. The graph of the minimal resolution of it is shown in Figure~\ref{fig2.2}.

On this surface we will make two examples for two different projections (pencils).
Firstly, let $\pi =(f,g): (Z,z)\to (\C^2,0)$ be defined by
$f(x,y,z)= u\circ
\pi= x$ and $g(x,y,z)= v\circ \pi=y$.
The discriminant curve of $\pi$ is the curve $v(u^2+v^4)=0$. Notice that this projection is not a
generic one because the image of the curve $\{ g= 0\}$ is  an irreducible component of the discriminant
curve and the image of $\{ f=0\}$ is tangent to the discriminant curve.

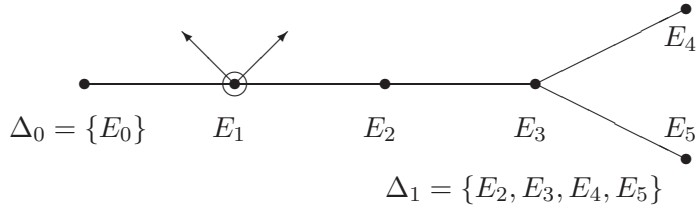
\begin{figure}[h]
$$\unitlength=1.00mm
\hspace*{-80mm}\begin{picture}(100.00,20.00)(0,10)
\thinlines

\put(50,15){\circle*{1.5}}
\put(50,15){\line(1,0){60}}
\put(70,15){\circle*{1.5}}
\put(70,15){\circle{3}}
\put(70,15){\vector(1,1){7}}
\put(70,15){\vector(-1,1){7}}
\put(90,15){\circle*{1.5}}
\put(110,15){\circle*{1.5}}
\put(130,25){\circle*{1.5}}
\put(130,5){\circle*{1.5}}
\put(110,15){\line(2,1){20}}
\put(110,15){\line(2,-1){20}}
\put(40,8){$\Delta_0=\{E_0\}$}
\put(67,8){$E_1$}
\put(87,8){$E_2$}
\put(107,8){$E_3$}
\put(127,20){$E_4$}
\put(127,8){$E_5$}
\put(90,0){$\Delta_1= \{ E_2,E_3,E_4,E_5\}$}
\end{picture}
$$
\bigskip
\caption{The graph of  the minimal good resolution of $\Lambda$.
} \label{fig2.3}
\end{figure}

The minimal good resolution of $\Lambda$ is just equal to the one of $(Z,z)$ and
there exists a unique dicritical component $E_1$: the divisor with weight
$(1,-2,0)$.
Thus, one has two special zones,
$SZ(\Lambda)=\{\Delta_0 ,\Delta_1\}$ with  $\Delta_0=\{ E_0\}$ and  $\Delta_1=\{ E_2, E_3, E_4,E_5\}$
(see Figure \ref{fig2.3} for the notations).
The Hironaka quotients corresponding to each vertex are:
$q(E_0)=q(E_1)=1$, $q(E_2)=3/2$ and $q(E_3)=q(E_4)=q(E_5)=2$.

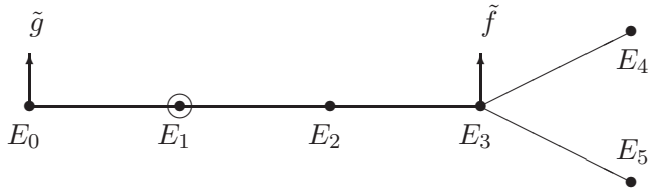
\begin{figure}[h]
\vspace*{10mm}
$$\unitlength=1.00mm
\hspace*{-80mm}\begin{picture}(100.00,20.00)(0,1)
\thinlines
\put(50,15){\circle*{1.5}}
\put(70,15){\circle{3}}
\put(50,15){\line(1,0){60}}
\put(70,15){\circle*{1.5}}
\put(50,15){\vector(0,1){7}}
\put(50,25){$\tilde g$}
\put(47,10){$E_0$}
\put(67,10){$E_1$}
\put(90,15){\circle*{1.5}}
\put(87,10){$E_2$}
\put(110,15){\circle*{1.5}}
\put(130,25){\circle*{1.5}}
\put(107,10){$E_3$}
\put(110,15){\vector(0,1){7}}
\put(110,25){$\tilde f$}
\put(130,5){\circle*{1.5}}
\put(128,20){$E_4$}
\put(128,8){$E_5$}
\put(110,15){\line(2,1){20}}
\put(110,15){\line(2,-1){20}}
\end{picture}
$$
\caption{The graph of  the minimal good resolution of $(f,g)$.
} \label{fig2.4}
\end{figure}

The connected component $\Delta _0$ doesn't contain any rupture component and $\Delta_1$ admits
a rupture component of Hironaka quotient equal to $2$.
The special fibre associated to $\Delta_1$ is $\{ f=0\}$ whose strict transform meets $\Delta _1$ at
$E_3$, and there are two irreducible components of $C(\pi)$ intersecting $\Delta_1$  at $E_4$ and
$E_5$.
The special fibre of $\Lambda$ associated to $\Delta _0$ is $\{g=0\}$ which is also a non reduced
irreducible component of the critical locus. It intersects $\Delta _0$ at $E_0$.
The minimal good resolution of the pencil $\Lambda$ is also the minimal good resolution of $(f,g)$, so
the corresponding graph
of the minimal good resolution of $fg=0$ is represented in figure \ref{fig2.4}.

\begin{figure}[h]
$$\unitlength=1.00mm
\hspace*{-80mm}\begin{picture}(100.00,20.00)(0,10)
\thinlines
\put(50,15){\circle*{1.5}}
\put(48,10){$E_0$}
\put(50,15){\line(1,0){60}}
\put(70,15){\circle*{1.5}}
\put(70,15){\circle{3}}
\put(70,15){\vector(0,1){7}}
\put(70,25){$\tilde f$}
\put(68,10){$E_1$}
\put(90,15){\circle*{1.5}}
\put(88,10){$E_2$}
\put(110,15){\circle*{1.5}}
\put(130,25){\circle*{1.5}}
\put(108,10){$E_3$}
\put(110,15){\vector(0,1){7}}
\put(110,25){$\tilde g$}
\put(130,5){\circle*{1.5}}
\put(128,20){$E_4$}
\put(128,8){$E_5$}
\put(110,15){\line(2,1){20}}
\put(110,15){\line(2,-1){20}}
\end{picture}
$$
\caption{The graph of  the minimal good resolution of $(f,g)$.
} \label{fig2.2bis}
\end{figure}
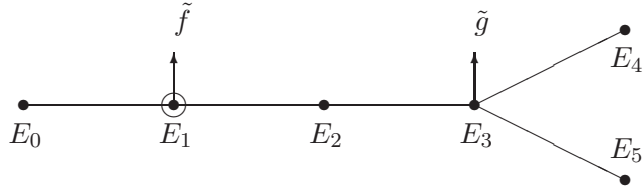

\bigskip
For the second example on $D_6$, let
the projection $\pi =(f,g): (Z,z)\to (\C^2,0)$  defined by $f(x,y,z)= x+y=u$ and $g(x,y,z)= x+2iy^2=v$.
As in the previous one the minimal good resolution of $\Lambda$ and the one of $\{fg=0\}$ coincides with
the minimal good resolution of $(Z,z)$. However, now
the graph of the minimal good resolution of $\pi$ is sligthly different and it is  represented in figure
\ref{fig2.2bis}.

In this case $f$ is a generic element of the pencil $\Lambda$ and $g$ is the special element associated
to $\Delta_1$. The special fibre of $\Lambda$ associated to $\Delta _0$ is  $g-f=0$. It is also a non reduced
irreducible component of the critical locus $C(\pi )$.
In this case the Hironaka quotients are
$q(E_0)=q(E_1)=1$, $q(E_2)=2/3$ and $q(E_3)=q(E_4)=q(E_5)=1/2$.

 \subsection{Example 3}

With this example, issued from \cite{LMW}, we illustrate the case where a special zone is a
singular point of the dicritical locus.

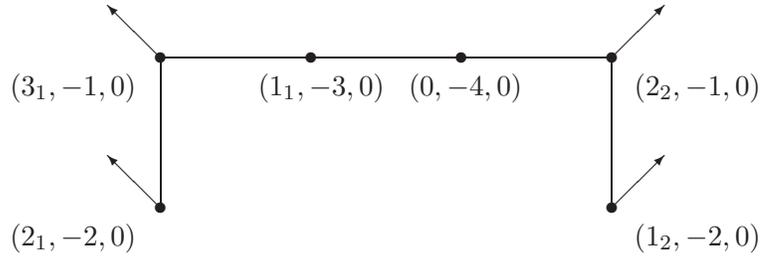
\begin{figure}[h]
\vspace*{10mm}
\unitlength=1.00mm
$$
\hspace*{-30mm}\begin{picture}(100.00,20.00)(0,10)
\thinlines
\put(50,35){\circle*{1.5}}
\put(50,15){\circle*{1.5}}
\put(30,30){$(3_1,-1,0)$}
\put(113,30){$(2_2,-1,0)$}
\put(63,30){$(1_1,-3,0)$}
\put(83,30){$(0,-4,0)$}
\put(30,10){$(2_1,-2,0)$}
\put(113,10){$(1_2,-2,0)$}
\put(50,35){\line(1,0){60}}
\put(50,35){\line(0,-1){20}}
\put(110,35){\line(0,-1){20}}
\put(70,35){\circle*{1.5}}
\put(90,35){\circle*{1.5}}
\put(110,35){\circle*{1.5}}
\put(110,15){\circle*{1.5}}
\put(50,35){\vector(-1,1){7}}
\put(110,35){\vector(1,1){7}}
\put(50,15){\vector(-1,1){7}}
\put(110,15){\vector(1,1){7}}
\end{picture}
$$
\caption{Graph of  the discriminant of $\pi$ and the coordinates axes.
} \label{fig3.1}
\end{figure}

 Let $(Z,z)$ be defined by
$z^2 =(x^2+y^5)(y^2+x^3)$ and let $\pi =(f,g): (Z,z)\to (\C^2,0)$ be the projection on the $(x,y)$-plane.
The  dual graph of the minimal embedded resolution of the discriminant curve  $(u^2+v^5)(v^2+u^3)=0$ of
$\pi$ and the coordinate axes   is shown in Figure \ref{fig3.1}.

\begin{figure}[h]
\vspace*{30mm}
\unitlength=1.0mm
$$
\hspace*{-30mm}\begin{picture}(100.00,20.00)(0,05)
\thinlines
\put(60,15){\circle*{1.5}}
\put(60,45){\circle*{1.5}}
\put(90,15){\circle*{1.5}}
\put(90,45){\circle*{1.5}}
\put(90,15){\circle{3}}
\put(90,45){\circle{3}}
\put(50,10){$(1^1,-2,0)$}
\put(80,10){$(0^1,-3,0)$}
\put(50,50){$(1^2,-2,0)$}
\put(75,50){$(0^2,-3,0)$}
\put(30,30){$\Delta_0=\{E_{1^1}, E_{1^2}\}$}
\put(95,30){$\Delta_1=E_{0^1}\cap E_{0^2}$}
\put(60,15){\line(1,0){40}}
\put(90,15){\line(0,1){30}}
\put(60,15){\line(0,1){30}}
\put(60,45){\line(1,0){40}}
\put(90,15){\vector(1,1){8}}
\put(90,15){\vector(1,0){10}}
\put(90,45){\vector(1,1){8}}
\put(90,45){\vector(1,0){10}}
\end{picture}
$$
\caption{The graph of  the minimal good resolution of $\Lambda$.
} \label{fig3.2}
\end{figure}
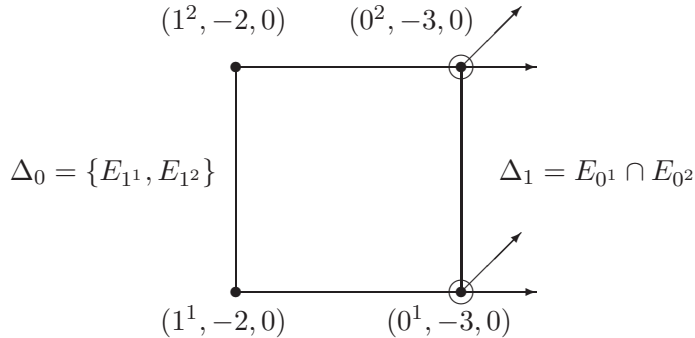

The graph of the minimal good resolution of $\Lambda$ is in figure \ref{fig3.2}.
The components $E_{0^1}$ and $E_{0^2}$ are dicritical. Thus, there exists two
special zones $\Delta_0$ and
$\Delta_1$ with $\Delta_0= \{ E_{1^1},E_{1^2}\}$ and $\Delta _1=
E_{0^1}\cap E_{0^2}= \{P\}$ where  $P$  is the singular point of $\cal D$.

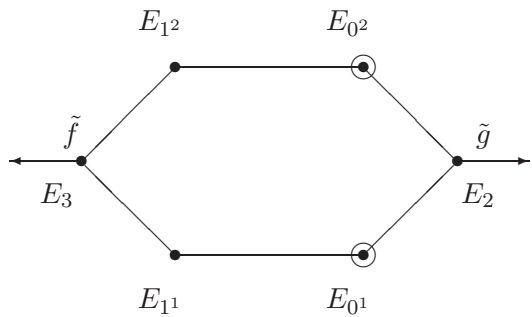
\begin{figure}[h]
\vspace*{20mm}
\unitlength=1.00mm
$$
\hspace*{-30mm}\begin{picture}(100.00,20.00)(0,10)
\thinlines
\put(50,15){\circle*{1.5}}
\put(50,40){\circle*{1.5}}
\put(75,15){\circle*{1.5}}
\put(75,40){\circle*{1.5}}
\put(75,15){\circle{3}}
\put(75,40){\circle{3}}
\put(37.5,27.5){\circle*{1.5}}
\put(87.5,27.5){\circle*{1.5}}
\put(45,8){$E_{1^1}$}
\put(70,8){$E_{0^1}$}
\put(45,45){$E_{1^2}$}
\put(70,45){$E_{0^2}$}
\put(88,22){$E_2$}
\put(32,22){$E_3$}
\put(50,15){\line(1,0){25}}
\put(75,15){\line(1,1){12}}
\put(75,40){\line(1,-1){12}}
\put(50,15){\line(-1,1){12}}
\put(50,40){\line(-1,-1){12}}
\put(50,40){\line(1,0){25}}
\put(87,27.5){\vector(1,0){10}}
\put(38,27.5){\vector(-1,0){10}}
\put(90,30){$\tilde g$}
\put(35,30){$\tilde f$}

\end{picture}
$$
\caption{The graph of  the minimal good resolution of $(f,g)$.
} \label{fig3.3}
\end{figure}

The special fibres associated to $\Delta_0$ and $\Delta_1$ are respectively $\{f=0\}$ and $\{ g=0\}$
and the graph of the minimal good resolution of $(f,g)$ is in figure \ref{fig3.3}.

The Hironaka quotients of the rational components (of self-intersection $-1$) $E_2$
and $E_3$ are respectively $2/3$ and $5/2$ and there exists two
irreducible components of $C(\pi)$ whose strict transform intersects $E_2$ and two
others whose strict transform intersects $E_3$.

\begin{thebibliography}{19}

\bibitem{B1} R. Bondil, {\it Discriminant of a generic projection of a minimal normal surface singularity}, C. R. Acad. Sci. Paris S\'er. I Math. {\bf 337} (2003), p. 195-200.

\bibitem{B2} R. Bondil, {\it General elements of an $m$-primary ideal on a normal surface singularity}, S\'eminaires et congr\`es {\bf 10} (2005), p. 11-20.

\bibitem{BL} R. Bondil and  L\^e D\~ung Tr\`ang, {\it Caract\'erisation des \'el\'ements superficiels d'un id\'eal}, C. R. Acad. Sci. Paris S\'er. I Math. {\bf 332} (2001), p. 717-722.

\bibitem{BNP} L. Birbrair, W. Neumann, A. Pichon, {\it The thick-thin decomposition and the bilipschitz classification of normal surface singularities}, preprint.


\bibitem{BPV}  W. Barth, C. Peters,  and A. Van de Ven, {\it Compact complex surfaces}, Ergebnisse der Mathematik, Springer-Verlag (1984).

\bibitem{Ca} E. Casas-Alvero, {\it Singularities of Plane Curves}, London Math. Soc. Lecture
Note Ser. 276, Cambridge Univ. Press, 2000.

\bibitem{DM}  F. Delgado and  H. Maugendre, {\it Special fibres and critical locus for a  pencil of plane curve
singularities}, Compositio Math. 136, 69--87 (2003).


\bibitem{Greuel_Buchweitz} R.-O Buchweitz and G.-M Greuel,  {\it The Milnor number and deformations of complex curve singularities}, Invent. Math. 58 (1980), 241-281.

\bibitem{L1} H. Laufer, Normal two dimensional singularities, Ann. of Math. Studies {\bf 71}, (1971), Princeton Univ. Press.

\bibitem{L2} H. Laufer,   {\it On normal two-dimensional double point singularities}, Israel Journal of Math., vol.31, n$^¡$ 3-4, 315-334 (1978).


\bibitem{LW1} L\^e D\~ung Tr\`ang and C. Weber, {\it \'Equisingularit\'e dans les pinceaux de germes de courbes planes et $C^0$-suffisance}, L'enseignement math\'ematique, {\bf 43}, (1997), 355-380.

\bibitem{LW2} L\^e D\~ung Tr\`ang and C. Weber,  {\it R\'esoudre est un jeu d'enfants}, Sem. Inst. de Estud. con Iberoamerica y Portugal, Tordesillas (1998).



\bibitem{LMW} L\^e D\~ung Tr\`ang, H. Maugendre and  C. Weber, {\it Geometry of critical loci}, Journal of the L.M.S. { \bf 63} (2001), 533-552.



\bibitem{Mi}  F. Michel,  {\it
Jacobian curves for normal complex surfaces},
Brasselet, J-P. (ed.) et al., Singularities II. Geometric and topological aspects. Proceedings of the international conference ``School and workshop on the geometry and topology of singularities" in honor of the 60th birthday of
L\^e D\~ung Tr\`ang, Cuernavaca, Mexico, January 8-26, 2007. Providence, RI: American Mathematical Society (AMS). Contemporary Mathematics 475, 135-150 (2008).

\bibitem{Mu} D. Mumford, {\it The topology of normal singularities of an algebraic surface and a criterion for simplicity}, Publications Math\'ematiques de l'IHES, tome {\bf 9} (1961), 5-22.

\bibitem{MM} H. Maugendre and F. Michel, {\it Fibrations associ\'ees \`a un pinceau de germes de courbes planes}, Ann. Fac. Sci. Toulouse, S\'er. 6, vol. X, fasc. 4, (2001), 745-777.

\bibitem{N} W. Neumann, {\it A calculus for plumbing applied to the topology of complex surface singularities and degenerating complex curves}, Trans. AMS {\bf 268}, (1981), 299-344.

\bibitem{S} J. Snoussi, {\it Limites d'espaces tangents \`a une surface normale}, Comment. Math. Helv. {\bf 73}, (2001), 61-88.

\bibitem{T} Teissier {\it The hunting of invariants on the geometry of discriminants}. Proc. of the Nordic
Summer School "Real and Complex Singularities", Oslo 1976, Sijthoff and Noordhooff 1977.

\bibitem{Wahl} J. Wahl, {\it Topology, geometry and equations of normal surface singularities},  in Singularities and Computer Algebra, LMS Lecture Notes Series {\bf 324}, Cam. Univ. Press (2006), 351-372.

\end{thebibliography}
\end{document}